\documentclass[10pt]{amsart}

\usepackage{amssymb,amsthm,amsmath}
\usepackage[numbers,sort&compress]{natbib}
\usepackage{color}
\usepackage{graphicx}
\usepackage{tikz}
\usepackage{comment}

\title[Stability of KAM tori ]{A note on stability of Eliasson-Kuksin's KAM tori for the
nonlinear Schr\"{o}dinger equation }

\author{Xiaolong He}
\address[Xiaolong He]
{Department of mathematics, Hangzhou Normal University, Hangzhou, 311121, China}
 \email{xlhe@hznu.edu.cn, hexiaolong@hnu.edu.cn}

\author{Jia Shi}
\address[Jia Shi]
{School of Mathematical Sciences,
Fudan University,
Shanghai 200433, China} \email{15110180007@fudan.edu.cn}

\author{Xiaoping Yuan}
\address[Xiaoping Yuan]
{School of Mathematical Sciences,
Fudan University,
Shanghai 200433, China} \email{xpyuan@fudan.edu.cn}


\theoremstyle{plain}
\newtheorem{thm}{Theorem}[section]
 
 \newtheorem{lem}[thm]{Lemma}
 \newtheorem{prop}[thm]{Proposition}
 \newtheorem{defn}[thm]{Definition}
 \newtheorem{rem}[thm]{Remark}
 
 \numberwithin{equation}{section}

\begin{document}


\begin{abstract}

Eliasson and Kuksin developed a KAM
approach to study the persistence of the invariant tori for  nonlinear Schr\"{o}dinger equation on $\mathbb{T}^{d}$. In this note, we
 improve Eliasson and Kuksin's KAM theorem   by using Kolmogorov's iterative scheme and  obtain a local normal form
for the transformed Hamiltonian.  As a consequence, we are able to derive the time $\delta^{-1}$ stability of the obtained
KAM tori.

\end{abstract}
\maketitle
\tableofcontents

\section{Introduction}\label{Sect 1}

In this paper, we consider the existence and the long time stability of the invariant tori for the $d$-dimensional nonlinear Schr\"odinger equation (NLS)
\begin{equation}\label{1}
-\mathrm{i}\dot{u}=-\Delta u+V(x)\ast u+\varepsilon\frac{\partial F}{\partial \bar{u}}(x,u,\bar{u}), \quad u=u(t,x)
\end{equation}
under the periodic boundary condition $x\in\mathbb{T}^{d}$.
The convolution function
$V: \mathbb{T}^{d}\rightarrow \mathbb{C}$ is analytic and the Fourier coefficient $\hat{V}(a)$
takes real value,
when expanding $V$ into Fourier series
$V(x)=\sum_{a\in\mathbb{Z}^{d}}\hat{V}(a)e^{\mathrm{i}\langle a,x\rangle}$.
The nonlinearity $F$ is real analytic in $x, \Re u, \Im u$.

The NLS equation \eqref{1} can be written as an infinite
dimensional Hamiltonian system. The KAM theory for Hamiltonian PDEs  started
in late 1980's and originally applied  to the  one dimensional PDEs,
which is now pretty well understood (see for example \cite{K,KP,P1,P2,W}).
However, the space-multidimensional KAM theory for Hamiltonian PDEs is
at its early stage. The first breakthrough in this respect was made by Bourgain
\cite{B1} on the two dimensional NLS equation, in which he developed further
Craig and Wayne's scheme on the persistency problem  of periodic solutions.
Later, Bourgain and his collaborators developed  new techniques of Green's function estimates
for linear problems, based on which he proved the persistence of invariant tori for space-multidimensional NLS and NLW equations \cite{B2}.
This method is now known as the Craig-Wayne-Bourgain (CWB) method.   However, Bourgain's original proof
gives no information on the linear stability of the obtained invariant tori. Recently, there are some progress
in applying the CWB method to obtain the linear stability of KAM tori for finitely dimensional Hamiltonian system \cite{HSSY},
but it remains an open question for the infinitely dimensional case.

The classical KAM approach for space-multidimensional Hamiltonian PDEs is developed by
Eliasson and Kuksin in their  paper \cite{EK} on NLS equation. They take  a sequence of symplectic transformations such that the transformed Hamiltonian guarantees the existence of the invariant torus.
Moreover, the KAM approach in \cite{EK} provides also the reducibility and linear stability of the obtained invariant
tori. However, they have to pay for the price that the number of the  second Melnikov conditions  becomes infinite when solving the block-diagonal homological equation,
which is far more complicated than that in the one-dimensional case.
To reduce those infinitely many conditions, they  analyze carefully the functions with the  T\"{o}plitz-Lipschitz property.
Moreover, due to the size of the blocks grows  much faster than quadratically along the iterations, they need to take sufficiently many
normal form computations at each KAM step to obtain a much faster iteration scheme.
See also \cite{EGK} and \cite{Y} for the KAM approach on the space-multidimensional  beam equation and  some
shallow-water equations, respectively.  But the problem on the linear stability of the KAM tori for space-multidimensional
nonlinear wave equation is still open and requires special attention.

Once the existence of the invariant tori of the Hamiltonian PDEs is established, one naturally concerns  the stability problem of the obtained KAM tori.
It is known that there exists transfer of energy in the cubic defocusing NLS equation \cite{CK}, which shows the instability phenomenon in Hamiltonian PDEs.
In the present paper, however, we concentrate on the long time stability of the KAM tori for the NLS equation \eqref{1}.

Consider a nonlinear differential equation
\begin{equation}\label{01}
  \dot{x}=X(x),
\end{equation}
which has an invariant tori $\mathcal{T}$ carrying the quasi-periodic flow $x(t)=x_{0}(t)$.
We say that the invariant tori $\mathcal{T}$ is linearly stable if the equilibrium of the linearized equation
\begin{equation*}
  \dot{y}=D X(x_{0}(t)) y
\end{equation*}
of \eqref{01} along $\mathcal{T}$
is Lyapunov stable.
However, in general, we cannot determine the stability of the nonlinear system \eqref{01}  from the  linear stability of the linearized equation.
This prompts us to study the nonlinear stability, among which the long time stability plays an important role in PDEs.

There are lots of literatures devoting to the long time stability of the equilibrium for the Hamiltonian PDEs.
Fortunately, it does not cause too much trouble to generalize the results in one dimensional to the multidimensional space (see for example \cite{BG}),
which is due to the fact that frequency shift does not come into small divisors during the normal form computations.
But not so for  the long time stability of KAM tori.
In \cite{CLY}, the long time stability of KAM tori for one dimensional NLS equation \eqref{1} (i.e., $d=1$) is proved by using the
tame property and the Birkhoff normal form. More precisely, given any $M>0$, there exists small $\delta_{0}>0$ such that for any $0<\delta<\delta_{0}$ and any solution $u(x,t)$ of  \eqref{1}
with its initial value $u(x,0)$ staying in a $\delta$-neighborhood (under the Sobolev norm) of some KAM tori $\mathcal{T}$ , the solution
$u(x,t)$ will stay inside the $C \delta$-neighborhood of $\mathcal{T}$ for all $0<|t|< \delta^{-M}$.

When it comes to $d\geq 2$ for  equation \eqref{1}, we are not able to apply the results in \cite{EK}  directly.
The reason is that the domain  (of the action variable and the normal coordinate) of the transformed Hamiltonian decreases to zero after the KAM iteration,
which does not affect the existence and linear stability  of invariant tori, but is not helpful to study the long time stability.
To this end, we improve Eliasson-Kuksin's KAM theorem, by employing Kolmogorov's iteration scheme (see Remark \ref{Kol} below)
and solving the modified homological equation as in \cite{CLY},
 to obtain a uniform domain for the transformed Hamiltonian,
which enables us to get a local normal form around the obtained KAM tori.
Moreover, based on the local normal form, it is easy to show the time $\delta^{-1}$ stability of KAM tori.
However, due to the much more complicated small divisor problems and the frequency shift, we  are not able to establish
the time $\delta^{-M}$ (with large $M$) stability in this note, although  the tame property
for space-multidimensional NLS equation can be preserved during the KAM iteration.

The main result of this paper is as follows.

\begin{thm}\label{t}
Under the assumptions for equation \eqref{1},
if $\varepsilon>0$ is sufficiently small, then for most $V$ (in the sense of measure), the $d$-dimensional nonlinear Schr\"odinger equation \eqref{1}
has a quasi-periodic solution. Moreover, assume $u_{0}(t,x)$ with initial value $u_{0}(0,x)$ is a quasi-periodic solution for the equation \eqref{1}.
Then, for any solution $u(t,x)$ with initial value $u(0,x)$ satisfying
\begin{equation*}
\|u(0,\cdot)-u_{0}(0,\cdot)\|_{H^{p}(\mathbb{T}^{d})}<\delta, \ \forall~ 0<\delta\ll1,
\end{equation*}
we have
\begin{equation*}
\|u(t,\cdot)-u_{0}(t,\cdot)\|_{H^{p}(\mathbb{T}^{d})}<C\delta, \ \forall~ 0<|t|<\delta^{-1}.
\end{equation*}
\end{thm}

In order to prove Theorem \ref{t}, we write the nonlinear Schr\"odinger equation (\ref{1}) as an infinite dimensional Hamiltonian system.
We keep the notations consistent with those in \cite{EK}.
Write
\begin{equation*}
u(x)=\sum_{a\in\mathbb{Z}^{d}}u_{a}e^{\mathrm{i}\langle a,x\rangle}, \ \overline{u(x)}=\sum_{a\in\mathbb{Z}^{d}}v_{a}e^{\mathrm{i}\langle -a,x\rangle},
\end{equation*}
and let
\begin{equation*}
\zeta_{a}=\left(
        \begin{array}{c}
          \xi_{a} \\
          \eta_{a} \\
        \end{array}
      \right)=\frac{1}{\sqrt{2}}\left(
                                  \begin{array}{c}
                                    u_{a}+v_{a} \\
                                    -\mathrm{i}(u_{a}-v_{a}) \\
                                  \end{array}
                                \right).
\end{equation*}
Then the nonlinear Schr\"odinger equation (\ref{1}) becomes a real Hamiltonian system
with the symplectic structure $d\xi\wedge d\eta$ and the Hamiltonian
\begin{equation*}
  \frac{1}{2}\sum_{a\in\mathbb{Z}^{d}}(|a|^{2}+\hat{V}(a))(\xi_{a}^{2}+\eta_{a}^{2})+\varepsilon \int_{\mathbb{T}^{d}}F(x,u(x), \overline{u(x)})dx.
\end{equation*}
Let $\mathcal{A}$ be a finite subset of $\mathbb{Z}^{d}$ and $\mathcal{L}=\mathbb{Z}^{d}\setminus\mathcal{A}$.
Introduce action angle variables $(\varphi_{a},r_{a})$, $a\in \mathcal{A}$,
\begin{equation*}
\xi_{a}=\sqrt{2(r_{a}+q_{a})}\cos\varphi_{a}, \ \eta_{a}=\sqrt{2(r_{a}+q_{a})}\sin\varphi_{a}, \ q_{a}>0.
\end{equation*}
Let
\begin{equation*}
\omega_{a}=|a|^{2}+\hat{V}(a), a\in \mathcal{A}, \ \Omega_{a}=|a|^{2}+\hat{V}(a), a\in \mathcal{L}.
\end{equation*}
We have the Hamiltonian
\begin{equation*}
h+f=\sum_{a\in \mathcal{A}}\omega_{a}r_{a}+\frac{1}{2}\sum_{a\in\mathcal{L}}\Omega_{a}(\xi_{a}^{2}+\eta_{a}^{2})+\varepsilon \int_{\mathbb{T}^{d}}F(x,u(x), \overline{u(x)})dx.
\end{equation*}
Assume $f$ is real analytic on
\begin{equation*}
D(\rho,\mu,\sigma)=\{(\varphi,r,\zeta)\in(\mathbb{C}/2\pi\mathbb{Z})^{\mathcal{A}}\times\mathbb{C}^{\mathcal{A}}\times l^{2}_{p}:|\Im\varphi|\leq\rho,|r|\leq\mu,\|\zeta\|_p\leq\sigma\},
\end{equation*}
where
\begin{equation*}
\|\zeta\|^{2}_{p}=\sum_{a\in\mathcal{L}}(|\xi^{2}_{a}|+|\eta^{2}_{a}|)\langle a\rangle^{2p}, \ \langle a\rangle=\max(|a|,1).
\end{equation*}

In the KAM iteration, we have symplectic maps
\begin{equation*}
\Phi_{j}:D(\rho_{j+1},\mu_{j+1},\sigma_{j+1})\rightarrow D(\rho_{j},\mu_{j},\sigma_{j})
\end{equation*}
such that $(h_{j}+f_{j})\circ\Phi_{j}=h_{j+1}+f_{j+1}$.
The existence of the KAM tori follows from the transformed Hamiltonian vector
field on the uniform domain $\cap_{j\geq 0} D(\rho_{j}, \mu_{j},\sigma_{j} )$.
In Eliasson-Kuksin \cite{EK}, there is $\rho_{j}\rightarrow \frac{\rho}{2}$, $\mu_{j}\rightarrow 0$, $\sigma_{j}\rightarrow 0$, which
causes some trouble to study the long time stability of the obtained KAM tori.
In this paper, we use Kolmogorov's iterative procedure as in \cite{CLY} such that
$\rho_{j}\rightarrow \frac{\rho}{2}$, $\mu_{j}\rightarrow \frac{\mu}{4}$, $\sigma_{j}\rightarrow \frac{\sigma}{2}$.
This is the main improvement on Eliasson and Kuksin's paper.
As a consequence, it is easy to show the time $\delta^{-1}$ stability of the obtained KAM tori.

The paper is organized as follows. In Section \ref{Sect 2}, we formulate and solve the homological equation in the KAM iteration.
Some definitions and notations are  introduced at the beginning.
In Section \ref{Sect 3}, following Eliasson-Kuksin \cite{EK}, we prove a KAM theorem for $d$-dimensional NLS,
based on which our main Theorem \ref{t} is an immediate result.

In this paper, $\|\cdot\|$ is an operator norm or $l^{2}$ norm. $|\cdot|$ will in general denote a sup norm.
For  $a \in \mathbb{Z}^{d}$, we use $|a|$ for the $l^{2}$ norm. The dimension $d$ will be fixed and $p>\frac{d}{2}$.
Let $\mathcal{A}$ be a finite subset of $\mathbb{Z}^{d}$ and $\mathcal{L}=\mathbb{Z}^{d}\setminus\mathcal{A}$.
Denote $\langle \zeta,\zeta'\rangle=\sum(\xi_{a}\xi'_{a}+\eta_{a}\eta'_{a})$ and $J=\left(
     \begin{array}{cc}
       0 & 1 \\
       -1 & 0 \\
     \end{array}
   \right)
$.

\section{Homological equation}\label{Sect 2}

In this section, we formulate and solve the homological equation in the KAM iteration.
To obtain an open and uniform domain for the transformed Hamiltonian, we apply Kolmogorov's iteration
scheme. As a result, the homological equation is complicated than that in \cite{EK}, but it can be solved
by the method developed in \cite{EK}. To begin with, we introduce some notations and definitions.

\subsection{Notations and definitions}
Recalling the Hamiltonian formulation of \eqref{1} in Section \ref{Sect 1}, we define   the
phase space and the associated norms.
For $\gamma\geq 0$, we denote
\begin{equation*}
l^{2}_{p,\gamma}=\{\zeta=(\xi,\eta)\in\mathbb{C}^{\mathcal{L}}\times\mathbb{C}^{\mathcal{L}}:\|\zeta\|_{p,\gamma}<\infty\},
\end{equation*}
where
\begin{equation*}
\|\zeta\|^{2}_{p,\gamma}=\sum_{a\in\mathcal{L}}(|\xi^{2}_{a}|+|\eta^{2}_{a}|)e^{2\gamma|a|}\langle a\rangle^{2p}, \ \langle a\rangle=\max(|a|,1).
\end{equation*}
When $\gamma= 0$, we write $l^{2}_{p}$ and $\|\zeta\|_p$ for simplicity.
The phase space of the Hamiltonian dynamical system is defined by
\begin{equation*}
\mathcal{P}^{p}=(\mathbb{C}/2\pi\mathbb{Z})^{\mathcal{A}}\times\mathbb{C}^{\mathcal{A}}\times l^{2}_{p}.
\end{equation*}
Let $U\subset[-1,1]^{\mathcal{A}}$ be a parameter set with positive measure.

\subsubsection{T\"{o}plitz-Lipschitz property}
We recall here some definitions and notations of the T\"{o}plitz-Lipschitz property in \cite{EK}. For more
heuristic formulation and explanation, we refer to \cite[Part I ]{EK}\cite{EKZ}.

Let $\mathrm{gl}(2,\mathbb{C})$ be the space of all complex $2\times2$-matrices.
For $A=\left(
         \begin{array}{cc}
           a & b \\
           c & d \\
         \end{array}
       \right)
\in \mathrm{gl}(2,\mathbb{C})$, denote $\pi A=\frac{1}{2}\left(
         \begin{array}{cc}
           a+d & b-c \\
           c-b & a+d \\
         \end{array}
       \right)
$ and $[A]=\left(
         \begin{array}{cc}
           |a| & |b| \\
           |c| & |d| \\
         \end{array}
       \right)
$.
Now consider an infinite-dimensional $\mathrm{gl}(2,\mathbb{C})$-valued matrix
\begin{equation*}
A:\mathcal{L}\times\mathcal{L}\rightarrow \mathrm{gl}(2,\mathbb{C}), \ (a,b)\mapsto A_{a}^{b}.
\end{equation*}
For reader's convenience, we collect the following notations and  definitions for the infinite-dimensional matrix $A$.
\begin{enumerate}
\item $|A|_\mathcal{D}$.  For any $\mathcal{D}\subset\mathcal{L}\times\mathcal{L}$, define
\begin{equation*}
|A|_\mathcal{D}=\sup_{(a,b)\in\mathcal{D}}\|A_{a}^{b}\|,
\end{equation*}
where $\|\cdot\|$ is the operator norm.

\item $|A|_{\gamma}$. Letting $(\pi A)_{a}^{b}=\pi A_{a}^{b}$ and $(\mathcal{E}_{\gamma}^{\pm}A)_{a}^{b}=[A_{a}^{b}]e^{\gamma|a\mp b|}$,
we define the norm
\begin{equation*}
|A|_\gamma=\max(|\mathcal{E}_{\gamma}^{+}\pi A|_{\mathcal{L}\times\mathcal{L}},\quad  |\mathcal{E}_{\gamma}^{-}(1-\pi) A|_{\mathcal{L}\times\mathcal{L}}).
\end{equation*}

\item $\mathcal{T}_{\Delta}^{\pm} A$ and $\mathcal{T}_{\Delta} A$.
We also define the truncation operators
\begin{equation*}
\mathcal{T}_{\Delta}^{\pm}A=A|_{\{(a,b)\in\mathcal{L}\times\mathcal{L}:|a\mp b|\leq\Delta\}}, \quad \mathcal{T}_{\Delta}A=\mathcal{T}_{\Delta}^{+}\pi A+\mathcal{T}_{\Delta}^{-}(1-\pi) A.
\end{equation*}

\item $A_{a}^{b}(\pm, c)$. A matrix $A:\mathcal{L}\times\mathcal{L}\rightarrow \mathrm{gl}(2,\mathbb{C})$ is \emph{T\"{o}plitz at $\infty$}, if for all $a,b,c$, the two limits
\begin{equation*}
\lim_{t\rightarrow +\infty}A_{a+tc}^{b\pm tc}\exists =A_{a}^{b}(\pm,c).
\end{equation*}

\item $D^{\pm}_{\Lambda}(c)$.
For $\Lambda\geq 0$, define the \emph{Lipschitz domain} $D_{\Lambda}^{+}(c)\subset\mathcal{L}\times\mathcal{L}$ to be the set of all $(a,b)$ such that
there exist $a',b'\in \mathbb{Z}^{d}, t\geq 0$ such that
\begin{equation*}
|a=a'+tc|\geq\Lambda(|a'|+|c|)|c|, \ |b=b'+tc|\geq\Lambda(|b'|+|c|)|c|, \ \frac{|a|}{|c|},\frac{|b|}{|c|}\geq 2\Lambda^{2}.
\end{equation*}
Define  $(a,b)\in D_{\Lambda}^{-}(c)$ if and only if $(a,-b)\in D_{\Lambda}^{+}(c)$.

\item $\mathrm{Lip}_{\Lambda,\gamma}^{\pm}A$.
For $c\neq 0$, we let $(\mathcal{M}_{c}A)_{a}^{b}=\left(\max(\frac{|a|}{|c|},\frac{|b|}{|c|})+1\right)[A_{a}^{b}]$ and
define the \emph{Lipschitz constants}
\begin{equation*}
\mathrm{Lip}_{\Lambda,\gamma}^{\pm}A=\sup_{c}|\mathcal{E}_{\gamma}^{\pm}\mathcal{M}_{c}(A-A(\pm,c))|_{D_{\Lambda}^{\pm}(c)}.
\end{equation*}

\item $\langle A\rangle_{\Lambda,\gamma}$. For $c\neq 0$,  the \emph{Lipschitz norm} is given by
\begin{equation*}
\langle A\rangle_{\Lambda,\gamma}=\max(\mathrm{Lip}_{\Lambda,\gamma}^{+}\pi A, ~\mathrm{Lip}_{\Lambda,\gamma}^{-}(1-\pi) A)+|A|_\gamma.
\end{equation*}

\item $|A|_{\gamma; U}$. Let $A(w):\mathcal{L}\times\mathcal{L}\rightarrow \mathrm{gl}(2,\mathbb{C})$ be $C^1$ (in the sense of Whitney) in $w\in U$, define
\begin{equation*}
|A|_{\gamma;U}=\sup_{w\in U}(|A(w)|_{\gamma},~ |\partial_{w}A(w)|_{\gamma}),\quad |A|_{U}=|A|_{0; U}.
\end{equation*}

\item $\langle A\rangle_{\Lambda, \gamma; U}$.  Let $A(w):\mathcal{L}\times\mathcal{L}\rightarrow \mathrm{gl}(2,\mathbb{C})$ be $C^1$ (in the sense of Whitney) in $w\in U$. If $A(w),\partial_{w}A(w)$ are T\"{o}plitz at $\infty$ for all $w\in U$, define
\begin{equation*}
\langle A\rangle_{\Lambda,\gamma;U}=\sup_{w\in U}(\langle A(w)\rangle_{\Lambda,\gamma},~\langle \partial_{w}A(w)\rangle_{\Lambda,\gamma}),\quad
\langle A\rangle_{\Lambda, U}=\langle A\rangle_{\Lambda,0; U}.
\end{equation*}

\item $d$-T\"{o}plitz. For $d=2$, the matrix $A$ is T\"{o}plitz-Lipschitz if it is T\"{o}plitz at $\infty$ and $\langle A\rangle_{\Lambda,\gamma}<\infty$ for some $\Lambda,\gamma$.
For $d>2$, we can define T\"{o}plitz-Lipschitz matrices inductively (see Section 2.4 in \cite{EK}).
\end{enumerate}

Based on  T\"{o}plitz-Lipschtiz matrices, we define the T\"{o}plitz-Lipschtiz property for functions.
\begin{defn}\label{def 2.1}
Let
$
D^{\gamma}(\sigma)=\{\zeta\in l^{2}_{p,\gamma}:\|\zeta\|_{p,\gamma}\leq\sigma\},
$
and the function $f:D^{0}(\sigma)\rightarrow \mathbb{C}$ be real analytic. We say
$f$ is T\"{o}plitz at $\infty$ if $\partial^{2}_{\zeta}f(\zeta)$ is T\"{o}plitz at $\infty$ for all $\zeta\in D^{0}(\sigma)$.
The norm $[f]_{\Lambda,\gamma,\sigma}$ is defined by  the smallest $C$ such that
\begin{equation*}
\begin{aligned}
& |f(\zeta)|\leq C,\quad \forall \zeta\in D^{0}(\sigma), \\
and \quad& \|\partial_{\zeta}f(\zeta)\|_{p,\gamma'}\leq\frac{C}{\sigma},\quad\langle\partial^{2}_{\zeta}f(\zeta)\rangle_{\Lambda,\gamma'}\leq\frac{C}{\sigma^{2}},
 \quad \forall \zeta\in D^{\gamma'}(\sigma), \quad \forall \gamma'\leq\gamma.
 \end{aligned}
\end{equation*}
\end{defn}
For more quantitative estimates of T\"{o}plitz-Lipschitz functions, we refer to \cite[Section 3]{EK}.


\subsubsection{The tame property}

Comparing with Eliassion-Kuksin \cite{EK}, we shall also
prove the preservation of the tame property of the Hamiltonian in the KAM iteration.
To begin with, we define the norms of analytic functions on the phase space.

\begin{defn}
Let
\begin{equation*}
D(\rho)=\{\varphi\in(\mathbb{C}/2\pi\mathbb{Z})^{\mathcal{A}}:|\Im\varphi|\leq\rho\},
\end{equation*}
$f:D(\rho)\times U\rightarrow \mathbb{C}$ be analytic in $\varphi\in D(\rho)$ and $C^1$  in $w\in U$,
\begin{equation*}
f(\varphi;w)=\sum_{k\in\mathbb{Z}^{\mathcal{A}}}\hat{f}(k;w)e^{\mathrm{i}\langle k,\varphi\rangle}.
\end{equation*}
Define the norm
\begin{equation*}
\|f\|_{D(\rho)\times U}=\sup_{w\in U}\sum_{k\in\mathbb{Z}^{\mathcal{A}}}\left(|\hat{f}(k;w)|+|\partial_{w}\hat{f}(k;w)|\right)e^{|k|\rho}.
\end{equation*}
\end{defn}

\begin{defn}
Let
\begin{equation}\label{D rho mu}
D(\rho,\mu)=\{(\varphi,r)\in(\mathbb{C}/2\pi\mathbb{Z})^{\mathcal{A}}\times\mathbb{C}^{\mathcal{A}}:|\Im\varphi|\leq\rho,|r|\leq\mu\},
\end{equation}
$f:D(\rho,\mu)\times U\rightarrow \mathbb{C}$ be analytic in $(\varphi,r)\in D(\rho,\mu)$ and $C^1$ in $w\in U$,
\begin{equation*}
f(\varphi,r;w)=\sum_{\alpha\in\mathbb{N}^{\mathcal{A}}}f^{\alpha}(\varphi;w)r^{\alpha}.
\end{equation*}
Define the norm
\begin{equation*}
\|f\|_{D(\rho,\mu)\times U}=\sum_{\alpha\in\mathbb{N}^{\mathcal{A}}}\|f^{\alpha}(\varphi;w)\|_{D(\rho)\times U}\mu^{|\alpha|}.
\end{equation*}
\end{defn}

\begin{defn}
Let
\begin{equation}\label{D rho mu sigma}
D(\rho,\mu,\sigma)=\{(\varphi,r,\zeta)\in(\mathbb{C}/2\pi\mathbb{Z})^{\mathcal{A}}\times\mathbb{C}^{\mathcal{A}}\times l^{2}_{p}:|\Im\varphi|\leq\rho,|r|\leq\mu,\|\zeta\|_p\leq\sigma\},
\end{equation}
$f:D(\rho,\mu,\sigma)\times U\rightarrow \mathbb{C}$ be analytic in $(\varphi,r,\zeta)\in D(\rho,\mu,\sigma)$ and $C^1$ in $w\in U$,
\begin{equation*}
f(\varphi,r,\zeta;w)=\sum_{\alpha\in\mathbb{N}^{\mathcal{A}},\beta\in\mathbb{N}^{\bar{\mathcal{L}}}}f^{\alpha\beta}(\varphi;w)r^{\alpha}\zeta^{\beta},
\end{equation*}
where $\bar{\mathcal{L}}=\mathcal{L}_{-1}\sqcup\mathcal{L},\mathcal{L}_{-1}=\mathcal{L}.$ For $a\in\mathcal{L}_{-1}, \zeta_{a}=\xi_{a}$, for $a\in\mathcal{L}, \zeta_{a}=\eta_{a}$.
Define the modulus
\begin{equation*}
\lfloor f \rceil_{D(\rho,\mu)\times U}(\zeta)=\sum_{\beta\in\mathbb{N}^{\bar{\mathcal{L}}}}\|f^{\beta}(\varphi,r;w)\|_{D(\rho,\mu)\times U}\zeta^{\beta},
\end{equation*}
where
\begin{equation*}
f^{\beta}(\varphi,r;w)=\sum_{\alpha\in\mathbb{N}^{\mathcal{A}}}f^{\alpha\beta}(\varphi;w)r^{\alpha}.
\end{equation*}
\end{defn}

For a homogeneous polynomial $f(\zeta)$ of degree $h>0$, it is associated with a symmetric $h$-linear form $\tilde{f}(\zeta^{(1)},\ldots,\zeta^{(h)})$ such that
$\tilde{f}(\zeta,\ldots,\zeta)=f(\zeta)$. For a monomial
\begin{equation*}
f(\zeta)=f^{\beta}\zeta^{\beta}=f^{\beta}\zeta_{j_{1}}\cdots\zeta_{j_{h}},
\end{equation*}
we define
\begin{equation*}
\tilde{f}(\zeta^{(1)},\ldots,\zeta^{(h)})=\widetilde{f^{\beta}\zeta^{\beta}}=\frac{1}{h!}\sum_{\tau_{h}}f^{\beta}\zeta_{j_{1}}^{(\tau_{h}(1))}\cdots\zeta_{j_{h}}^{(\tau_{h}(h))},
\end{equation*}
where $\tau_{h}$ is an $h$-permutation. For a homogeneous polynomial
\begin{equation*}
f(\zeta)=\sum_{|\beta|=h}f^{\beta}\zeta^{\beta},
\end{equation*}
we have
\begin{equation*}
\tilde{f}(\zeta^{(1)},\ldots,\zeta^{(h)})=\sum_{|\beta|=h}\widetilde{f^{\beta}\zeta^{\beta}}.
\end{equation*}

Now we can define \textbf{$p$-tame norm} of a Hamiltonian vector field. We first consider a Hamiltonian
\begin{equation}\label{f_h}
f(\varphi,r,\zeta;w)=f_{h}=\sum_{\alpha\in\mathbb{N}^{\mathcal{A}},\beta\in\mathbb{N}^{\bar{\mathcal{L}}},|\beta|=h}f_{h}^{\alpha\beta}(\varphi;w)r^{\alpha}\zeta^{\beta},
\end{equation}
which is homogenous polynomial in $\zeta$ of order $h$.
Letting $f_{\zeta}=(f_{\eta},-f_{\xi})$, the Hamiltonian vector field $X_f$ is given by $(f_{r},-f_{\varphi},f_{\zeta})$.
For $h\geq 1$, denote
\begin{equation*}
\|(\zeta^h)\|_{p,1}=\frac{1}{h}\sum_{j=1}^{h}\|\zeta^{(1)}\|_{1}\cdots\|\zeta^{(j-1)}\|_{1}\|\zeta^{(j)}\|_{p}\|\zeta^{(j+1)}\|_{1}\cdots\|\zeta^{(h)}\|_{1}.
\end{equation*}

Recall the domains $D(\rho,\mu)$ and $D(\rho,\mu,\sigma)$ in \eqref{D rho mu} and \eqref{D rho mu sigma},
respectively.
\begin{defn}
For $f$ in \eqref{f_h}, we define
\begin{equation*}
|||f_{\zeta}|||^{T}_{p,D(\rho,\mu)\times U}=\sup_{0\neq(\zeta^{j})\in l^2_p, 1\leq j\leq h-1}\frac{\|\widetilde{\lfloor f_{\zeta} \rceil}_{D(\rho,\mu)\times U}
(\zeta^{(1)},\ldots,\zeta^{(h-1)})\|_{p}}{\|(\zeta^{h-1})\|_{p,1}}, \ h\geq 2,
\end{equation*}
\begin{equation*}
|||f_{\zeta}|||^{T}_{p,D(\rho,\mu)\times U}=\|\widetilde{\lfloor f_{\zeta} \rceil}_{D(\rho,\mu)\times U}\|_{p}, \ h=0,1.
\end{equation*}
Then we define the $p$-tame norm of $f_{\zeta}$ by
\begin{equation*}
|||f_{\zeta}|||^{T}_{p,D(\rho,\mu,\sigma)\times U}=\max(|||f_{\zeta}|||^{T}_{p,D(\rho,\mu)\times U}, ~|||f_{\zeta}|||^{T}_{1,D(\rho,\mu)\times U})~ \sigma^{h-1}.
\end{equation*}
\end{defn}

\begin{defn}
For $f$ in \eqref{f_h}, we define
\begin{equation*}
|||f_{r}|||_{D(\rho,\mu)\times U}=\sup_{0\neq(\zeta^{j})\in l^2_1, 1\leq j\leq h}\frac{|\widetilde{\lfloor f_{r} \rceil}_{D(\rho,\mu)\times U}
(\zeta^{(1)},\ldots,\zeta^{(h)})|}{\|(\zeta^{h})\|_{1,1}}, \ h\geq 1,
\end{equation*}
\begin{equation*}
|||f_{r}|||_{D(\rho,\mu)\times U}=|\widetilde{\lfloor f_{r} \rceil}_{D(\rho,\mu)\times U}|, \ h=0.
\end{equation*}
Then we define the  norm of $f_{r}$ by
\begin{equation*}
|||f_{r}|||_{D(\rho,\mu,\sigma)\times U}=|||f_{r}|||_{D(\rho,\mu)\times U}\sigma^{h}.
\end{equation*}
The  norm of $f_{\varphi}$ is defined as that of $f_{r}$.
\end{defn}

\begin{defn}
For $f$ in \eqref{f_h},
the $p$-tame norm of the Hamiltonian vector field $X_f$ is defined by
\begin{equation*}
|||X_f|||^{T}_{p,D(\rho,\mu,\sigma)\times U}=|||f_{r}|||_{D(\rho,\mu,\sigma)\times U}+\frac{1}{\mu}|||f_{\varphi}|||_{D(\rho,\mu,\sigma)\times U}+\frac{1}{\sigma}|||f_{\zeta}|||^{T}_{p,D(\rho,\mu,\sigma)\times U}.
\end{equation*}
\end{defn}

\begin{defn}
For a general Hamiltonian
\begin{equation*}
f(\varphi,r,\zeta;w)=\sum_{h\geq 0}f_{h}, \quad f_{h}=\sum_{\alpha\in\mathbb{N}^{\mathcal{A}},\beta\in\mathbb{N}^{\bar{\mathcal{L}}},|\beta|=h}f_{h}^{\alpha\beta}(\varphi;w)r^{\alpha}\zeta^{\beta},
\end{equation*}
define the $p$-tame norm of the Hamiltonian vector field $X_f$ by
\begin{equation*}
|||X_f|||^{T}_{p,D(\rho,\mu,\sigma)\times U}=\sum_{h\geq 0}|||X_{f_{h}}|||^{T}_{p,D(\rho,\mu,\sigma)\times U}.
\end{equation*}
\end{defn}

\begin{rem}
The $p$-tame norm can also be defined in complex coordinates
\begin{equation*}
z=\left(
    \begin{array}{c}
      u \\
      v \\
    \end{array}
  \right)
  =C^{-1}\left(
           \begin{array}{c}
             \xi \\
             \eta \\
           \end{array}
         \right),
         C=\frac{1}{\sqrt{2}}\left(
         \begin{array}{cc}
           1 & 1 \\
           -\mathrm{i} & \mathrm{i} \\
         \end{array}
       \right).
\end{equation*}
\end{rem}

Following the proof of Theorem 3.1 in \cite{CLY}, we have the following proposition.

\begin{prop}\label{p}
If $0<\tau<\rho, 0<\tau'<\frac{\sigma}{2}$, then
\begin{equation*}
|||X_{\{f,g\}}|||^{T}_{p,D(\rho-\tau,(\sigma-\tau')^{2},\sigma-\tau')\times U}\leq
C\max\left(\frac{1}{\tau},\frac{\sigma}{\tau'}\right)|||X_{f}|||^{T}_{p,D(\rho,\sigma^{2},\sigma)\times U}|||X_{g}|||^{T}_{p,D(\rho,\sigma^{2},\sigma)\times U},
\end{equation*}
where $C>0$ is a constant depending on $\# \mathcal{A}$.
\end{prop}

We define the weighted norm of the Hamiltonian vector field $X_f$ by
\begin{equation*}
|||X_f|||_{\mathcal{P}^{p},D(\rho,\mu,\sigma)\times U}=\sup_{(\varphi,r,\zeta;w)\in D(\rho,\mu,\sigma)\times U}\|X_{f}\|_{\mathcal{P}^{p},D(\rho,\mu,\sigma)},
\end{equation*}
where
\begin{equation*}
\|X_{f}\|_{\mathcal{P}^{p},D(\rho,\mu,\sigma)}=|f_{r}|+\frac{1}{\mu}|f_{\varphi}|+\frac{1}{\sigma}\|f_{\zeta}\|_{p}.
\end{equation*}
Following the proof of Theorem 3.5 in \cite{CLY}, we have
\begin{equation*}
|||X_f|||_{\mathcal{P}^{p},D(\rho,\mu,\sigma)\times U}\leq|||X_f|||^{T}_{p,D(\rho,\mu,\sigma)\times U}.
\end{equation*}

\subsubsection{Normal form matrices}
We introduce the blocks in \cite{EK}.
For $\Delta\geq 0$, define an equivalence relation on  $\mathcal{L}=\mathbb{Z}^{d}\setminus \mathcal{A}$ generated by the pre-equivalence relation
\begin{equation*}
a\sim b \Leftrightarrow |a|=|b|,|a-b|\leq\Delta.
\end{equation*}
Let $[a]_\Delta$ be the equivalence class (block) of $a$ and $\mathcal{E}_\Delta$ be the set of equivalence classes.
Let $d_\Delta$ be the supremum of all block diameters, then by Proposition 4.1 in \cite{EK}, there is $d_\Delta\preceq \Delta^{\frac{(d+1)!}{2}}$.

A matrix $A:\mathcal{L}\times\mathcal{L}\rightarrow \mathrm{gl}(2,\mathbb{C})$ is on normal form, denoted $\mathcal{NF}_{\Delta}$, if
$A$ is real valued, symmetric, $\pi A=A$ and block-diagonal over $\mathcal{E}_\Delta$, i.e., $ A_{a}^{b}=0, \forall [a]_\Delta\neq [b]_\Delta$.
A matrix $Q:\mathcal{L}\times\mathcal{L}\rightarrow \mathbb{C}$ is on normal form, denoted $\mathcal{NF}_{\Delta}$, if
$Q$ is Hermitian and block-diagonal over $\mathcal{E}_\Delta$.

For a normal form matrix $A$, we write
\begin{equation*}
\frac{1}{2}\langle\zeta,A\zeta\rangle=\frac{1}{2}\langle\xi,A_{1}\xi\rangle+\langle\xi,A_{2}\eta\rangle+\frac{1}{2}\langle\eta,A_{1}\eta\rangle,
\end{equation*}
where $A_{1}+\mathrm{i}A_{2}$ is a Hermitian matrix.
Let
\begin{equation*}
z=\left(
    \begin{array}{c}
      u \\
      v \\
    \end{array}
  \right)
  =C^{-1}\left(
           \begin{array}{c}
             \xi \\
             \eta \\
           \end{array}
         \right),
         C=\frac{1}{\sqrt{2}}\left(
         \begin{array}{cc}
           1 & 1 \\
           -\mathrm{i} & \mathrm{i} \\
         \end{array}
       \right)
\end{equation*}
and define $C^{T}AC:\mathcal{L}\times\mathcal{L}\rightarrow \mathrm{gl}(2,\mathbb{C})$ by $(C^{T}AC)_{a}^{b}=C^{T}A_{a}^{b}C$.
If $A$ is on normal form, then
\begin{equation*}
\frac{1}{2}\langle z,C^{T}ACz\rangle=\langle u,Qv\rangle,
\end{equation*}
where $Q$ is the normal form matrix associated to $A$.

\subsection{The homological equations}
In this part, we formulate and solve the homological equations in the KAM iteration, from which
we derive the quantitive estimates of the symplectic transformation and the new Hamiltonian vector field.
See Proposition \ref{p1}.
  To begin with, we make some basic
assumptions \textbf{(H1)-(H3)} below.
\smallskip

\textbf{Assumption (H1):}
Consider the unperturbed Hamiltonian
\begin{equation*}
h(r,\zeta;w)=\langle\omega(w),r\rangle+\frac{1}{2}\langle\zeta,(\Omega(w)+H(w))\zeta\rangle,
\end{equation*}
where $\Omega(w)$ is a real diagonal matrix with diagonal elements $\Omega_{a}(w)I$,
$H(w),\partial_{w}H(w)$ are T\"{o}plitz at $\infty$ and $\mathcal{NF}_{\Delta}$ for all $w\in U$.

Assume
\begin{equation}\label{as1}
\partial_{w_{a}}\omega_{b}(w)=\delta_{ab}, \ a\in \mathcal{A}, b\in \mathcal{A}, w\in U,
\end{equation}
\begin{equation}\label{as2}
\partial_{w_{a}}\Omega_{b}(w)=\delta_{ab}, \ a\in \mathcal{A}, b\in \mathcal{L}, w\in U.
\end{equation}
Assume further that there exist constants $c_1,c_2,c_3,c_4,c_5>0$ such that
\begin{equation}\label{as3}
|\Omega_{a}(w)-|a|^{2}|\leq c_1e^{-c_2|a|}, \ a\in \mathcal{L}, w\in U,
\end{equation}
\begin{equation}\label{as4}
|\Omega_{a}(w)|\geq c_3, \ a\in \mathcal{L}, w\in U,
\end{equation}
\begin{equation}\label{as5}
|\Omega_{a}(w)+\Omega_{b}(w)|\geq c_3, \ a,b\in \mathcal{L}, w\in U,
\end{equation}
\begin{equation}\label{as6}
|\Omega_{a}(w)-\Omega_{b}(w)|\geq c_3, \   |a|\neq |b|,      \ a,b\in \mathcal{L}, w\in U,
\end{equation}
\begin{equation}\label{as7}
\|H(w)\|\leq \frac{c_3}{4}, \  w\in U,
\end{equation}
\begin{equation}\label{as8}
\|\partial_{w}H(w)\|\leq c_4, w\in U,
\end{equation}
\begin{equation}\label{as9}
\langle H \rangle_{\Lambda;U}\leq c_5.
\end{equation}

\textbf{Assumption (H2):}
Consider the perturbation $f: D(\rho, \mu, \sigma)\times U\rightarrow \mathbb{C}$ and write
\begin{equation*}
f(\varphi,r,\zeta;w)=f^{low}+f^{high},
\end{equation*}
where
\begin{equation}\label{f low}
f^{low}= f^{\varphi}+f^{0}+f^{1}+f^{2}
= F^{\varphi}(\varphi;w)+\langle F_{0}(\varphi;w),r\rangle+\langle F_{1}(\varphi;w),\zeta\rangle+\frac{1}{2}\langle F_{2}(\varphi;w)\zeta,\zeta\rangle.
\end{equation}
Assume that
\begin{equation}\label{ho1}
|||X_{f^{low}}|||^{T}_{p,D(\rho,\mu,\sigma)\times U}\leq\varepsilon, \quad |||X_{f^{high}}|||^{T}_{p,D(\rho,\mu,\sigma)\times U}\leq1,
\end{equation}
Let
\begin{equation*}
D^{\gamma}(\rho,\mu,\sigma)=\{(\varphi,r,\zeta)\in(\mathbb{C}/2\pi\mathbb{Z})^{\mathcal{A}}\times\mathbb{C}^{\mathcal{A}}\times l^{2}_{p,\gamma}:|\Im\varphi|\leq\rho,|r|\leq\mu,\|\zeta\|_{p,\gamma}\leq\sigma\},
\end{equation*}
and suppose  that
$f:D^{\gamma}(\rho,\mu,\sigma)\times U\rightarrow \mathbb{C}$ is real analytic in $(\varphi,r,\zeta)\in D^{\gamma}(\rho,\mu,\sigma)$ and $C^1$ in $w\in U$.
Define
\begin{equation*}
[f]_{\Lambda,\gamma,\sigma;U,\rho,\mu} =\sup_{(\varphi,r)\in D(\rho,\mu)}[f(\varphi,r,\cdot;\cdot)]_{\Lambda,\gamma,\sigma;U},
\end{equation*}
in which the norm inside the supremum is defined in Definition \ref{def 2.1}.
Assume  further that
\begin{equation}\label{ho2}
[f^{low}]_{\Lambda,\gamma,\sigma;U,\rho,\mu}\leq\varepsilon, \quad [f^{high}]_{\Lambda,\gamma,\sigma;U,\rho,\mu}\leq 1.
\end{equation}
\smallskip

\textbf{Assumption (H3):} Let $\Delta'>1$ and $0<\kappa<1$. Assume there exists $U' \subset U$ such that
for all $w\in U'$, $0<|k|\leq\Delta'$,
\begin{equation}\label{sd1}
|\langle k,\omega(w)\rangle|\geq \kappa,
\end{equation}
\begin{equation}\label{sd2}
|\langle k,\omega(w)\rangle+\alpha(w)|\geq \kappa, \quad \forall~\alpha(w)\in\mathrm{spec}(((\Omega+H)(w))_{[a]_\Delta}), \ \forall~ [a]_\Delta,
\end{equation}
\begin{equation}\label{sd3}
|\langle k,\omega(w)\rangle+\alpha(w)+\beta(w)|\geq \kappa, \quad \forall
\left\{
\begin{aligned}
\alpha(w)\in\mathrm{spec}(((\Omega+H)(w))_{[a]_\Delta}),\\
\beta(w)\in\mathrm{spec}(((\Omega+H)(w))_{[b]_\Delta}),
\end{aligned}
\right.
~\forall~ [a]_\Delta, [b]_\Delta,
\end{equation}
 and
\begin{equation}\label{sd4}
|\langle k,\omega(w)\rangle+\alpha(w)-\beta(w)|\geq \kappa, \quad
 \forall
 \left\{\begin{aligned}
 \alpha(w)\in\mathrm{spec}(((\Omega+H)(w))_{[a]_\Delta}),\\
 \beta(w)\in\mathrm{spec}(((\Omega+H)(w))_{[b]_\Delta}),
 \end{aligned}\right.
\end{equation}
for $\mathrm{dist}([a]_\Delta, [b]_\Delta)\leq \Delta'+2d_\Delta.$

Recalling $f^{low}$ in \eqref{f low},
we define the truncation operator
\begin{equation*}
\mathcal{T}_{\Delta'}f^{low}=
\sum_{|k|\leq\Delta'}\left(\hat{F^{\varphi}}(k;w)+\langle\hat{F_{0}}(k;w),r\rangle+\langle\hat{F_{1}}(k;w),\zeta\rangle
+\frac{1}{2}\langle\mathcal{T}_{\Delta'}\hat{F_{2}}(k;w)\zeta,\zeta\rangle\right)e^{\mathrm{i}\langle k,\varphi\rangle}.
\end{equation*}
Now we state our main results on the homological equations, whose proof is delayed to the next subsection.

\begin{prop}\label{p1}

Under the assumptions \emph{\textbf{(H1)-(H3)}}, if
 $\gamma,\sigma,\rho,\mu<1$, $\Lambda,\Delta\geq 3$, $\rho=\sigma$, $\mu=\sigma^{2}$, $d_{\Delta}\gamma\leq1$,
then for all $w\in U'$, the homological equation
\begin{equation}\label{ho3}
\{h,s\}=-\mathcal{T}_{\Delta'}f^{low}-\mathcal{T}_{\Delta'}\{f^{high},s\}^{low}+h_1
\end{equation}
has solutions
\begin{equation}\label{ho4}
s(\varphi,r,\zeta;w)=s^{low}=s^{\varphi}+s^{0}+s^{1}+s^{2},
\end{equation}
\begin{equation}\label{ho5}
h_{1}(r,\zeta;w)=a_{1}(w)+\langle \chi_{1}(w),r\rangle+\frac{1}{2}\langle\zeta,H_{1}(w)\zeta\rangle
\end{equation}
with the estimates
\begin{equation}\label{ho6}
\begin{aligned}
& |||X_{s^{\varphi}}|||^{T}_{p,D(\rho-\tau,\sigma^{2},\sigma)\times U'}\preceq \frac{\varepsilon}{\tau\kappa^{2}},\\
& |||X_{s^{1}}|||^{T}_{p,D(\rho-3\tau,(\sigma-3\tau)^{2},\sigma-3\tau)\times U'}
\preceq\frac{d_{\Delta}^{d}\varepsilon}{\tau^{3}\kappa^{4}},\\
& |||X_{s^{0}}|||^{T}_{p,D(\rho-5\tau,(\sigma-5\tau)^{2},\sigma-5\tau)\times U'}\preceq \frac{d_{\Delta}^{d}\varepsilon}{\tau^{5}\kappa^{6}},\\
& |||X_{s^{2}}|||^{T}_{p,D(\rho-5\tau,(\sigma-5\tau)^{2},\sigma-5\tau)\times U'}\preceq
\frac{d_{\Delta}^{3d}\varepsilon}{\tau^{5}\kappa^{6}},\\
& |||X_{s}|||^{T}_{p,D(\rho-5\tau,(\sigma-5\tau)^{2},\sigma-5\tau)\times U'}\preceq\frac{d_{\Delta}^{3d}\varepsilon}{\tau^{5}\kappa^{6}},\\
& |||X_{h_{1}}|||^{T}_{p,D(\rho-5\tau,(\sigma-5\tau)^{2},\sigma-5\tau)\times U'}
\preceq\frac{d_{\Delta}^{d}\varepsilon}{\tau^{4}\kappa^{4}},
\end{aligned}
\end{equation}
where $0<\tau<\frac{\rho}{100}$, $a\preceq b$ means there exists a constant $c>0$ depending on
$d,\# \mathcal{A}, p, c_1,c_2,c_3,c_4,c_5$ such that $a\leq cb$.

The new Hamiltonian
\begin{equation}\label{ho12}
(h+f)\circ X^{t}_{s}|_{t=1}=h+h_1+f_1,
\end{equation}
\begin{equation}\label{ho13}
f_1=(1-\mathcal{T}_{\Delta'})f^{low}+f^{high}+(1-\mathcal{T}_{\Delta'})\{f^{high},s\}^{low}+\{f^{high},s\}^{high}
\end{equation}
\begin{equation*}
+\int_{0}^{1}(1-t)\{\{h,s\},s\}\circ X^{t}_{s}dt+\int_{0}^{1}\{f^{low},s\}\circ X^{t}_{s}dt+\int_{0}^{1}(1-t)\{\{f^{high},s\},s\}\circ X^{t}_{s}dt
\end{equation*}
with the estimates
\begin{equation}\label{ho14}
|||X_{f_{1}^{low}}|||^{T}_{p,D(\rho-8\tau,(\sigma-8\tau)^{2},\sigma-8\tau)\times U'}
\preceq\frac{d_{\Delta}^{d}e^{-\frac{1}{2}\tau\Delta'}}{\kappa^{4}\tau^{\# \mathcal{A}+4}}\varepsilon
+\frac{(\Delta\Delta')^{\exp}}{\sigma^{6}\kappa^{4}}\frac{e^{-\frac{1}{2}\gamma\Delta'}}{\gamma^{d+p}\tau^{\# \mathcal{A}+3}}\varepsilon+\frac{d_{\Delta}^{6d}\varepsilon^{2}}{\tau^{12}\kappa^{12}},
\end{equation}
\begin{equation}\label{ho15}
|||X_{f_{1}^{high}}|||^{T}_{p,D(\rho-8\tau,(\sigma-8\tau)^{2},\sigma-8\tau)\times U'}\preceq 1+\frac{d_{\Delta}^{3d}\varepsilon}{\tau^{6}\kappa^{6}}+\frac{d_{\Delta}^{6d}\varepsilon^{2}}{\tau^{12}\kappa^{12}},
\end{equation}
where the exponent $\exp$ depends on $d,\# \mathcal{A}, p$.

Moreover, the following estimates hold:

\emph{(1)} The functions $s$ and $h_{1}$ satisfy
\begin{equation}\label{ho16}
\begin{aligned}
& [s]_{\Lambda'+d_\Delta+2,\gamma,\sigma';U',\rho',\mu'}
\preceq\frac{1}{\kappa^{7}}(\Delta\Delta')^{\exp}\frac{1}{\rho-\rho'}\left(\frac{1}{\sigma-\sigma'}\frac{1}{\sigma}+\frac{1}{\rho-\rho'}\frac{1}{\mu-\mu'}\right)
\frac{1}{\mu}\varepsilon,\\
& [h_{1}]_{\Lambda'+d_\Delta+2,\gamma,\sigma';U',\rho',\mu'}
\preceq\frac{1}{\kappa^{4}}(\Delta\Delta')^{\exp}\frac{1}{\rho-\rho'}\left(\frac{1}{\sigma-\sigma'}\frac{1}{\sigma}+\frac{1}{\rho-\rho'}
\frac{1}{\mu-\mu'}\right)\frac{1}{\mu}\varepsilon,\\
\end{aligned}
\end{equation}
where $\rho'<\rho$, $\mu'<\mu$, $\sigma'<\sigma$, $\Lambda'\geq \mathrm{cte}.\max(\Lambda,d^{2}_{\Delta},d^{2}_{\Delta'})$,
the constant $\mathrm{cte}$. is the one in Proposition 6.7 in \cite{EK}.

\emph{(2)} There is measure estimate
\begin{equation}\label{ho18}
\mathrm{meas}(U\backslash U')\preceq\max(\Lambda,\Delta,\Delta')^{\exp}\kappa^{(\frac{1}{d+1})^{d}}.
\end{equation}

\end{prop}

\begin{rem}\label{Kol}
In Eliasson-Kuksin \cite{EK}, the homological equation therein is
\begin{equation*}
\{h, s\}=-\mathcal{T}_{\Delta'} f^{low}+h_{1}.
\end{equation*}
As we see, the term $\{f^{high}, s\}^{low}$ is not eliminated and is put into the new perturbation along the iteration. To obtain a fast convergent iteration,
they requires the radius of domain for the variables $r$ and $\zeta$ to decrease to zero. However, in the present paper, we follow Kolmogorov's iterative
scheme, which eliminates also the lower order term $\{f^{high}, s\}^{low}$ \emph{(}up to the truncation\emph{)}.
This certainly complicates the homological equation when expressing in Fourier modes \emph{(}see \eqref{es4}-\eqref{es7} below\emph{)},
but enables us to obtain a uniform domain after the KAM iteration. Based on the uniform domain, we can establish a
normal form of the Hamiltonian vector field around the obtained KAM tori and then study its local stability.
\end{rem}

\subsection{Proof of Proposition \ref{p1}}
The whole subsection is devoted to the proof of Proposition \ref{p1}.
Note  that the measure estimate (\ref{ho18}) follows from Proposition 6.6 and 6.7 in \cite{EK}, which
is not repeated here. Since the proof is too long, we divide it into five steps.
In what follows, we shall not indicate the dependence on the parameter $w$ of functions when it is known from the text.

\subsubsection{\textbf{Step 1: Write the homological equation in Fourier modes}}

Recall that
\begin{equation*}
f^{low}=f^{\varphi}+f^{0}+f^{1}+f^{2}=F^{\varphi}(\varphi;w)+\langle F_{0}(\varphi;w),r\rangle+\langle F_{1}(\varphi;w),\zeta\rangle+\frac{1}{2}\langle F_{2}(\varphi;w)\zeta,\zeta\rangle,
\end{equation*}
\begin{equation*}
s=s^{low}=s^{\varphi}+s^{0}+s^{1}+s^{2}=S^{\varphi}(\varphi;w)+\langle S_{0}(\varphi;w),r\rangle+\langle S_{1}(\varphi;w),\zeta\rangle+\frac{1}{2}\langle S_{2}(\varphi;w)\zeta,\zeta\rangle.
\end{equation*}
By the calculation of Section 4.1.2 in \cite{CLY}, we obtain
\begin{equation}
\label{es1}
\begin{aligned}
&\{f^{high},s\}^{low}=\{f^{high},s\}^{0}+\{f^{high},s\}^{1}+\{f^{high},s\}^{2},\\
&\{f^{high},s\}^{1}=\{f^{high},s^{\varphi}\}^{1},\\
&\{f^{high},s\}^{0}=\{f^{high},s^{\varphi}+s^{1}\}^{0},   \ \{f^{high},s\}^{2}=\{f^{high},s^{\varphi}+s^{1}\}^{2}.
\end{aligned}
\end{equation}
Let $g=\{f^{high},s\}$ and write
\begin{equation*}
g^{low}=g^{0}+g^{1}+g^{2}=\langle G_{0}(\varphi;w),r\rangle+\langle G_{1}(\varphi;w),\zeta\rangle+\frac{1}{2}\langle G_{2}(\varphi;w)\zeta,\zeta\rangle.
\end{equation*}

In Fourier modes, the homological equation (\ref{ho3}) decomposes into
\begin{equation}\label{es4}
-\mathrm{i}\langle k,\omega(w)\rangle\hat{S}^{\varphi}(k;w)=-\hat{F}^{\varphi}(k;w)+\delta_{0}^{k}a_{1}(w),
\end{equation}
\begin{equation}\label{es5}
-\mathrm{i}\langle k,\omega(w)\rangle\hat{S}_{1}(k;w)+J(\Omega(w)+H(w))\hat{S}_{1}(k;w)=-\hat{F}_{1}(k;w)-\hat{G}_{1}(k;w),
\end{equation}
\begin{equation}\label{es6}
-\mathrm{i}\langle k,\omega(w)\rangle\hat{S}_{0}(k;w)=-\hat{F}_{0}(k;w)-\hat{G}_{0}(k;w)+\delta_{0}^{k}\chi_{1}(w),
\end{equation}
\begin{equation}\label{es7}
-\mathrm{i}\langle k,\omega(w)\rangle\hat{S}_{2}(k;w)+(\Omega(w)+H(w))J\hat{S}_{2}(k;w)-\hat{S}_{2}(k;w)J(\Omega(w)+H(w))
\end{equation}
\begin{equation*}
=-\hat{F}_{2}(k;w)-\hat{G}_{2}(k;w)+\delta_{0}^{k}H_{1}(w).
\end{equation*}

\subsubsection{\textbf{Step 2:
Solve the homological equations \emph{(\ref{es4})-(\ref{es7}).}}}
Since equations \eqref{es4}-\eqref{es7} are coupled, we solve them in the following order
 $$(\ref{es4}) \rightarrow (\ref{es5}) \rightarrow (\ref{es6}) \rightarrow (\ref{es7}).$$

\textbf{\underline{Solution of \eqref{es4}.}} The homological equation \eqref{es4} is very standard in KAM theory. We obtain from (\ref{es4}) that
\begin{equation*}
a_{1}=\hat{F}^{\varphi}(0),~\textrm{and}~\hat{S}^{\varphi}(k)=\frac{\hat{F}^{\varphi}(k)}{\mathrm{i}\langle k,\omega\rangle},~k\neq0.
\end{equation*}
By (\ref{sd1}), we have
\begin{equation}\label{es9}
|\hat{S}^{\varphi}(k)|\leq\frac{1}{\kappa}|\hat{F}^{\varphi}(k)|.
\end{equation}
Differentiating (\ref{es4}) with respect to the parameter $w$, we derive a similar homological equation
\begin{equation}
\label{es10}
-\mathrm{i}\partial_{w}(\langle k,\omega(w)\rangle)\hat{S}^{\varphi}(k;w)-\mathrm{i}\langle k,\omega(w)\rangle\partial_{w}\hat{S}^{\varphi}(k;w)=-\partial_{w}\hat{F}^{\varphi}(k;w)
\end{equation}
for $\partial_{w} \hat{S}^{\varphi}(k;w)$
and there is also
\begin{equation*}
|\partial_{w}\hat{S}^{\varphi}(k)|\leq\frac{1}{\kappa}(|k||\hat{S}^{\varphi}(k)|+|\partial_{w}\hat{F}^{\varphi}(k)|),
\end{equation*}
which together with \eqref{es9} implies
\begin{equation*}
|\hat{S}^{\varphi}(k)|+|\partial_{w}\hat{S}^{\varphi}(k)|\leq\frac{|k|+1}{\kappa^{2}}(|\hat{F}^{\varphi}(k)|+|\partial_{w}\hat{F}^{\varphi}(k)|).
\end{equation*}
It then follows that
\begin{equation*}
\begin{aligned}
\|S^{\varphi}_{\varphi}&\|_{D(\rho-\tau)\times U'}=\sum_{k\in\mathbb{Z}^{\mathcal{A}}}(|\hat{S}^{\varphi}_{\varphi}(k)|+|\partial_{w}\hat{S}^{\varphi}_{\varphi}(k)|)e^{|k|(\rho-\tau)}\\
 &\leq\sum_{k\in\mathbb{Z}^{\mathcal{A}}}\frac{|k|+1}{\kappa^{2}}(|\hat{F}^{\varphi}_{\varphi}(k)|+|\partial_{w}\hat{F}^{\varphi}_{\varphi}(k)|)e^{|k|(\rho-\tau)}
 \preceq \frac{1}{\tau\kappa^{2}}\|F^{\varphi}_{\varphi}\|_{D(\rho)\times U'}.
\end{aligned}
\end{equation*}
Consequently, we have
\begin{equation}
\label{es14}
|||X_{s^{\varphi}}|||^{T}_{p,D(\rho-\tau,\sigma^{2},\sigma)\times U'}\preceq \frac{1}{\tau\kappa^{2}}|||X_{f^{\varphi}}|||^{T}_{p,D(\rho,\sigma^{2},\sigma)\times U'}.
\end{equation}
\smallskip

\textbf{\underline{Solution of \eqref{es5}.}} Now we consider the equation (\ref{es5}) and write it simply, we as
\begin{equation*}
\mathrm{i}\langle k,\omega(w)\rangle S+J(\Omega+H)S=F+G.
\end{equation*}
We change to complex coordinates
\begin{equation*}
z=\left(
    \begin{array}{c}
      u \\
      v \\
    \end{array}
  \right)
  =C^{-1}\left(
           \begin{array}{c}
             \xi \\
             \eta \\
           \end{array}
         \right),
         C=\frac{1}{\sqrt{2}}\left(
         \begin{array}{cc}
           1 & 1 \\
           -\mathrm{i} & \mathrm{i} \\
         \end{array}
       \right).
\end{equation*}
Let $S'=C^{-1}S=\left(
           \begin{array}{c}
             S'_{1} \\
             S'_{2} \\
           \end{array}
         \right),F'=C^{-1}F=\left(
           \begin{array}{c}
             F'_{1} \\
             F'_{2} \\
           \end{array}
         \right),G'=C^{-1}G=\left(
           \begin{array}{c}
             G'_{1} \\
             G'_{2} \\
           \end{array}
         \right)$.
The equation \eqref{es5}
becomes
\begin{equation}\label{es16}
\begin{aligned}
& \mathrm{i}\langle k,\omega(w)\rangle S'_{1}-\mathrm{i}(\Omega+H^{T})S'_{1}=F'_{1}+G'_{1},\\
& \mathrm{i}\langle k,\omega(w)\rangle S'_{2}+\mathrm{i}(\Omega+H)S'_{2}=F'_{2}+G'_{2}.
\end{aligned}
\end{equation}
We only solve $S_{1}'$ in \eqref{es16} since $S_{2}'$ can be solved accordingly.

By (\ref{sd2}), we have
\begin{equation*}
\|S'_{1[a]_\Delta}\|\leq\frac{1}{\kappa}\|F'_{1[a]_\Delta}+G'_{1[a]_\Delta}\|.
\end{equation*}
Using similar arguments to that of \eqref{es10}, we get
\begin{equation}\label{es21}
\|S'_{1[a]_\Delta}\|+\|\partial_{w}S'_{1[a]_\Delta}\|\preceq\frac{|k|+1}{\kappa^{2}}(\|F'_{1[a]_\Delta}+G'_{1[a]_\Delta}\|+\|\partial_{w}F'_{1[a]_\Delta}+\partial_{w}G'_{1[a]_\Delta}\|).
\end{equation}
For $z=(u,v)$, define $\tilde{z}=(\tilde{u},\tilde{v})$ such that for all $a\in [a]_\Delta$, $\tilde{u}_{a}=\|u_{[a]_\Delta}\|,\tilde{v}_{a}=\|v_{[a]_\Delta}\|$.
It follows from  (\ref{es21}) that
\begin{equation*}
\begin{aligned}
& |\sum_{a\in\mathcal{L}}\|S'_{1a\varphi}(\varphi)\|_{D(\rho-3\tau)\times U'}u_{a}| \\
\leq &
\sum_{a\in\mathcal{L}}\sum_{k\in\mathbb{Z}^{\mathcal{A}}}(|\hat{S}'_{1a\varphi}(k)|+|\partial_{w}\hat{S}'_{1a\varphi}(k)|)e^{(\rho-3\tau)|k|}|u_{a}|
\\
=&\sum_{k\in\mathbb{Z}^{\mathcal{A}}}\sum_{[a]_\Delta\in\mathcal{E}_{\Delta}}\sum_{a\in[a]_\Delta}(|\hat{S}'_{1a\varphi}(k)|+|\partial_{w}\hat{S}'_{1a\varphi}(k)|)|u_{a}|e^{(\rho-3\tau)|k|}  \\
\preceq&\sum_{k\in\mathbb{Z}^{\mathcal{A}}}\sum_{a\in\mathcal{L}}\frac{|k|+1}{\kappa^{2}}
(|(\hat{F}'_{1a\varphi}+\hat{G}'_{1a\varphi})(k)|+|(\partial_{w}\hat{F}'_{1a\varphi}+\partial_{w}\hat{G}'_{1a\varphi})(k) |)\tilde{u}_{a}e^{(\rho-3\tau)|k|}  \\
\preceq&\frac{1}{\tau\kappa^{2}}\sum_{a\in\mathcal{L}}\|F'_{1a\varphi}(\varphi)+G'_{1a\varphi}(\varphi)\|_{D(\rho-2\tau)\times U'}\tilde{u}_{a} .
\end{aligned}
\end{equation*}
There is a similar estimate for that of $\sum_{a\in\mathcal{L}}\|S'_{2a\varphi}(\varphi)\|_{D(\rho-3\tau)\times U'} v_{a}$.

Noticing that
\begin{equation*}
s^{1}=\sum_{a\in\mathcal{L}}(S'_{1a}(\varphi)u_{a}+S'_{2a}(\varphi)v_{a}),
\end{equation*}
and
\begin{equation*}
\lfloor s^{1}_{\varphi}\rceil_{D(\rho-3\tau,\sigma^{2})\times U'}(z)=\sum_{a\in\mathcal{L}}(\|S'_{1a\varphi}(\varphi)\|_{D(\rho-3\tau)\times U'}u_{a}+\|S'_{2a\varphi}(\varphi)\|_{D(\rho-3\tau)\times U'}v_{a}),
\end{equation*}
%
%
we get
\begin{equation*}
|\lfloor s^{1}_{\varphi}\rceil_{D(\rho-3\tau,\sigma^{2})\times U'}(z)|\preceq\frac{1}{\tau\kappa^{2}}\lfloor f^{1}_{\varphi}+g^{1}_{\varphi}\rceil_{D(\rho-2\tau,\sigma^{2})\times U'}(\tilde{z}).
\end{equation*}
It then follows from
\begin{equation*}
|||s^{1}_{\varphi}|||_{D(\rho-3\tau,\sigma^{2})\times U'}=\sup_{0\neq z\in l_{1}^{2}}\frac{|\lfloor s^{1}_{\varphi}\rceil_{D(\rho-3\tau,\sigma^{2})\times U'}(z)|}{\|z\|_{1}}
\end{equation*}
and $\|\tilde{z}\|_{1}\preceq d_{\Delta}^{d}\|z\|_{1}$ that
\begin{equation}\label{es28}
|||s^{1}_{\varphi}|||_{D(\rho-3\tau,\sigma^{2})\times U'}\preceq\frac{d_{\Delta}^{d}}{\tau\kappa^{2}}|||f^{1}_{\varphi}+g^{1}_{\varphi}|||_{D(\rho-2\tau,\sigma^{2})\times U'}.
\end{equation}

Next we estimate
\begin{equation}\label{es29}
\begin{aligned}
 |||s^{1}_{z}|||^{T}_{p,D(\rho-3\tau,\sigma^{2})\times U'}=&\|\lfloor s^{1}_{z}\rceil_{D(\rho-3\tau,\sigma^{2})\times U'}\|_{p},\\
 =&\left(\sum_{a\in\mathcal{L}}(\|S'_{1a}(\varphi)\|^{2}_{D(\rho-3\tau)\times U'}+\|S'_{2a}(\varphi)\|^{2}_{D(\rho-3\tau)\times U'})\langle a\rangle^{2p}\right)^{1/2}.
\end{aligned}
\end{equation}
By (\ref{es21}), we see from the Minkowski's inequality that
\begin{equation*}
\begin{aligned}
 (\sum_{a\in[a]_\Delta}\|S'_{1a}(\varphi)&\|^{2}_{D(\rho-3\tau)\times U'})^{\frac{1}{2}}
=
[\sum_{a\in[a]_\Delta}(\sum_{k\in\mathbb{Z}^{\mathcal{A}}}(|\hat{S}'_{1a}(k)|+|\partial_{w}\hat{S}'_{1a}(k)|)e^{(\rho-3\tau)|k|})^{2}]^{\frac{1}{2}}\\
&\leq\sum_{k\in\mathbb{Z}^{\mathcal{A}}}[\sum_{a\in[a]_\Delta}(|\hat{S}'_{1a}(k)|+|\partial_{w}\hat{S}'_{1a}(k)|)^{2}]^{\frac{1}{2}}e^{(\rho-3\tau)|k|}\\
&\preceq\sum_{k\in\mathbb{Z}^{\mathcal{A}}}\sum_{a\in[a]_\Delta}\frac{|k|+1}{\kappa^{2}}
(|(\hat{F}'_{1a}+\hat{G}'_{1a})(k)|+|(\partial_{w}\hat{F}'_{1a}+\partial_{w}\hat{G}'_{1a})(k)|)e^{(\rho-3\tau)|k|}\\
&\preceq\frac{1}{\tau\kappa^{2}}\sum_{a\in[a]_\Delta}\|F'_{1a}(\varphi)+G'_{1a}(\varphi)\|_{D(\rho-2\tau)\times U'}.
\end{aligned}
\end{equation*}
As a result, there is
\begin{equation*}
\sum_{a\in[a]_\Delta}\|S'_{1a}(\varphi)\|^{2}_{D(\rho-3\tau)\times U'}\langle a\rangle^{2p}
\preceq\left(\frac{d_{\Delta}^{d}}{\tau\kappa^{2}}\right)^{2}\sum_{a\in[a]_\Delta}\|F'_{1a}(\varphi)+G'_{1a}(\varphi)\|^{2}_{D(\rho-2\tau)\times U'}\langle a\rangle^{2p}.
\end{equation*}
and similar estimate of $\sum_{a\in[a]_\Delta}\|S'_{2a}(\varphi)\|^{2}_{D(\rho-3\tau)\times U'}\langle a\rangle^{2p}$ holds.
Then we have
\begin{equation*}
|||s^{1}_{z}|||^{T}_{p,D(\rho-3\tau,\sigma^{2})\times U'}\preceq\frac{d_{\Delta}^{d}}{\tau\kappa^{2}}|||f^{1}_{z}+g^{1}_{z}|||^{T}_{p,D(\rho-2\tau,\sigma^{2})\times U'},
\end{equation*}
which together with (\ref{es28}) implies
\begin{equation}
\label{es35}
|||X_{s^{1}}|||^{T}_{p,D(\rho-3\tau,\sigma^{2},\sigma)\times U'}\preceq \frac{d_{\Delta}^{d}}{\tau\kappa^{2}}|||X_{f^{1}+g^{1}}|||^{T}_{p,D(\rho-2\tau,\sigma^{2},\sigma)\times U'}.
\end{equation}

\underline{\textbf{Solution of \eqref{es6}}.}
Solving  equation (\ref{es6}) as  equation (\ref{es4}), we obtain
\begin{equation*}
|||X_{s^{0}}|||^{T}_{p,D(\rho-5\tau,\sigma^{2},\sigma)\times U'}\preceq \frac{1}{\tau\kappa^{2}}|||X_{f^{0}+g^{0}}|||^{T}_{p,D(\rho-4\tau,\sigma^{2},\sigma)\times U'}.
\end{equation*}

\underline{\textbf{Solution of \eqref{es7}}.}
For simplicity, we write (\ref{es7}) as
\begin{equation*}
\mathrm{i}\langle k,\omega(w)\rangle S+(\Omega+H)JS-SJ(\Omega+H)=F+G-H_{1}.
\end{equation*}
We change to complex coordinates $z=C^{-1}\zeta$.
Let $S'=C^{T}SC=\left(
                  \begin{array}{cc}
                    S_{1}' &  S_{2}' \\
                    S'^{T}_{2} & S_{3}' \\
                  \end{array}
                \right)$
and $F'=C^{T}FC,G'=C^{T}GC$, $H_{1}'=C^{T}H_{1}C$. Then equation \eqref{es7}
 becomes
\begin{equation}\label{es38}
\mathrm{i}\langle k,\omega(w)\rangle S_{1}'+\mathrm{i}(\Omega+H)S_{1}'+\mathrm{i}S_{1}'(\Omega+H^{T})=F'_{1}+G'_{1},
\end{equation}
\begin{equation}\label{es39}
\mathrm{i}\langle k,\omega(w)\rangle S_{2}'+\mathrm{i}(\Omega+H)S_{2}'-\mathrm{i}S_{2}'(\Omega+H)=F'_{2}+G'_{2}-H_{12}',
\end{equation}
\begin{equation}\label{es40}
\mathrm{i}\langle k,\omega(w)\rangle S'^{T}_{2}-\mathrm{i}(\Omega+H^{T})S'^{T}_{2}+\mathrm{i}S'^{T}_{2}(\Omega+H^{T})=F'^{T}_{2}+G'^{T}_{2}-H'^{T}_{12},
\end{equation}
\begin{equation}\label{es41}
\mathrm{i}\langle k,\omega(w)\rangle S_{3}'-\mathrm{i}(\Omega+H^{T})S_{3}'-\mathrm{i}S_{3}'(\Omega+H)=F'_{3}+G'_{3}.
\end{equation}

For $k\neq0,H_{12}'=0$. By (\ref{sd4}), we have
\begin{equation*}
\|S'^{[b]_\Delta}_{2[a]_\Delta}\|\leq\frac{1}{\kappa}\|F'^{[b]_\Delta}_{2[a]_\Delta}+G'^{[b]_\Delta}_{2[a]_\Delta}\|.
\end{equation*}
Using similar argument to that of \eqref{es10}, we get
\begin{equation}\label{es45}
\|S'^{[b]_\Delta}_{2[a]_\Delta}\|+\|\partial_{w}S'^{[b]_\Delta}_{2[a]_\Delta}\|
\preceq\frac{|k|+1}{\kappa^{2}}(\|F'^{[b]_\Delta}_{2[a]_\Delta}+G'^{[b]_\Delta}_{2[a]_\Delta}\|
+\|\partial_{w}F'^{[b]_\Delta}_{2[a]_\Delta}+\partial_{w}G'^{[b]_\Delta}_{2[a]_\Delta}\|).
\end{equation}
For $k=0,H_{12}'=(F'_{2}+G'_{2})|_{\{(a,b)\in\mathcal{L}\times\mathcal{L}:|a|=|b|\}}$
and estimate (\ref{es45}) still holds.

Using (\ref{sd3}), we also  have
\begin{equation}\label{es46}
\begin{aligned}
& \|S'^{[b]_\Delta}_{1[a]_\Delta}\|+\|\partial_{w}S'^{[b]_\Delta}_{1[a]_\Delta}\|
\preceq\frac{|k|+1}{\kappa^{2}}(\|F'^{[b]_\Delta}_{1[a]_\Delta}+G'^{[b]_\Delta}_{1[a]_\Delta}\|
+\|\partial_{w}F'^{[b]_\Delta}_{1[a]_\Delta}+\partial_{w}G'^{[b]_\Delta}_{1[a]_\Delta}\|),\\
& \|S'^{[b]_\Delta}_{3[a]_\Delta}\|+\|\partial_{w}S'^{[b]_\Delta}_{3[a]_\Delta}\|
\preceq\frac{|k|+1}{\kappa^{2}}(\|F'^{[b]_\Delta}_{3[a]_\Delta}+G'^{[b]_\Delta}_{3[a]_\Delta}\|
+\|\partial_{w}F'^{[b]_\Delta}_{3[a]_\Delta}+\partial_{w}G'^{[b]_\Delta}_{3[a]_\Delta}\|).
\end{aligned}
\end{equation}

Recall that
\begin{equation*}\label{es48}
s^{2}=\frac{1}{2}\sum_{a,b\in\mathcal{L}}(S'^{b}_{1a}(\varphi)u_{a}u_{b}+2S'^{b}_{2a}(\varphi)u_{a}v_{b}+S'^{b}_{3a}(\varphi)v_{a}v_{b}).
\end{equation*}
We consider first $s^{2}_{\varphi}$
\begin{equation*}\label{es49}
\begin{aligned}
\lfloor s^{2}_{\varphi}\rceil_{D(\rho-5\tau,\sigma^{2})\times U'}(z)=&\frac{1}{2}
\sum_{a,b\in\mathcal{L}}(\|S'^{b}_{1a\varphi}(\varphi)\|_{D(\rho-5\tau)\times U'}u_{a}u_{b}\\
&+2\|S'^{b}_{2a\varphi}(\varphi)\|_{D(\rho-5\tau)\times U'}u_{a}v_{b}+\|S'^{b}_{3a\varphi}(\varphi)\|_{D(\rho-5\tau)\times U'}v_{a}v_{b}),
\end{aligned}
\end{equation*}
and  its associated multilinear form $\widetilde{s^{2}_{\varphi}}$
\begin{equation}\label{es50}
\begin{aligned}
&\widetilde{\lfloor s^{2}_{\varphi}\rceil}_{D(\rho-5\tau,\sigma^{2})\times U'}(z^{(1)},z^{(2)})=\frac{1}{2}
\sum_{a,b\in\mathcal{L}}(\|S'^{b}_{1a\varphi}(\varphi)\|_{D(\rho-5\tau)\times U'}u^{(1)}_{a}u^{(2)}_{b}\\
&\quad+\|S'^{b}_{2a\varphi}(\varphi)\|_{D(\rho-5\tau)\times U'}(u^{(1)}_{a}v^{(2)}_{b}+u^{(2)}_{a}v^{(1)}_{b})+\|S'^{b}_{3a\varphi}(\varphi)\|_{D(\rho-5\tau)\times U'}v^{(1)}_{a}v^{(2)}_{b}).
\end{aligned}
\end{equation}
By (\ref{es46}), we know that $|\sum_{a\in[a]_\Delta,b\in[b]_\Delta}\|S'^{b}_{1a\varphi}(\varphi)\|_{D(\rho-5\tau)\times U'}u^{(1)}_{a}u^{(2)}_{b}|$ is less than
\begin{align*}
&\quad\sum_{a\in[a],b\in[b]}\sum_{k\in\mathbb{Z}^{\mathcal{A}}}(|\hat{S}'^{b}_{1a\varphi}(k)|+|\partial_{w}\hat{S}'^{b}_{1a\varphi}(k)|)e^{(\rho-5\tau)|k|}|u^{(1)}_{a}u^{(2)}_{b}| \\
&\preceq\sum_{k\in\mathbb{Z}^{\mathcal{A}}}\frac{|k|+1}{\kappa^{2}}(\|\hat{F}'^{[b]}_{1\varphi[a]}(k)+\hat{G}'^{[b]}_{1\varphi[a]}(k)\|
+\|\partial_{w}\hat{F}'^{[b]}_{1\varphi[a]}(k)+\partial_{w}\hat{G}'^{[b]}_{1\varphi[a]}(k)\|)\|u^{(1)}_{[a]}\|\|u^{(2)}_{[b]}\|e^{(\rho-5\tau)|k|}\\
&\preceq\sum_{k\in\mathbb{Z}^{\mathcal{A}}}\sum_{a\in[a],b\in[b]}\frac{|k|+1}{\kappa^{2}}(|\hat{F}'^{b}_{1a\varphi}(k)+\hat{G}'^{b}_{1a\varphi}(k)|
+|\partial_{w}\hat{F}'^{b}_{1a\varphi}(k)+\partial_{w}\hat{G}'^{b}_{1a\varphi}(k)|)\tilde{u}^{(1)}_{a}\tilde{u}^{(2)}_{b}e^{(\rho-5\tau)|k|}\\
&\preceq\frac{1}{\tau\kappa^{2}}\sum_{a\in[a],b\in[b]}\|F'^{b}_{1a\varphi}(\varphi)+G'^{b}_{1a\varphi}(\varphi)\|_{D(\rho-4\tau)\times U'}\tilde{u}^{(1)}_{a}\tilde{u}^{(2)}_{b},
\end{align*}
which implies
\begin{equation*}\label{es52}
|\sum_{a,b\in\mathcal{L}}\|S'^{b}_{1a\varphi}(\varphi)\|_{D(\rho-5\tau)\times U'}u^{(1)}_{a}u^{(2)}_{b}|
\preceq\frac{1}{\tau\kappa^{2}}\sum_{a,b\in\mathcal{L}}\|F'^{b}_{1a\varphi}(\varphi)+G'^{b}_{1a\varphi}(\varphi)\|_{D(\rho-4\tau)\times U'}\tilde{u}^{(1)}_{a}\tilde{u}^{(2)}_{b}.
\end{equation*}
There are similar estimates  for the other two summation in the R.H.S of \eqref{es50},
Then we obtain
\begin{equation*}\label{es56}
|\widetilde{\lfloor s^{2}_{\varphi}\rceil}_{D(\rho-5\tau,\sigma^{2})\times U'}(z^{(1)},z^{(2)})|\preceq\frac{1}{\tau\kappa^{2}}
\widetilde{\lfloor f^{2}_{\varphi}+g^{2}_{\varphi}\rceil}_{D(\rho-4\tau,\sigma^{2})\times U'}(\tilde{z}^{(1)},\tilde{z}^{(2)}).
\end{equation*}
Since
\begin{equation*}\label{es57}
|||s^{2}_{\varphi}|||_{D(\rho-5\tau,\sigma^{2})\times U'}
=\sup_{0\neq z^{(1)},z^{(2)}\in l_{1}^{2}}\frac{|\widetilde{\lfloor s^{2}_{\varphi}\rceil}_{D(\rho-5\tau,\sigma^{2})\times U'}(z^{(1)},z^{(2)})|}{\|z^{(1)}\|_{1}\|z^{(2)}\|_{1}},
\end{equation*}
we see from $\|\tilde{z}\|_{1}\preceq d_{\Delta}^{d}\|z\|_{1}$ that
\begin{equation}\label{es58}
|||s^{2}_{\varphi}|||_{D(\rho-5\tau,\sigma^{2})\times U'}\preceq\frac{d_{\Delta}^{2d}}{\tau\kappa^{2}}|||f^{2}_{\varphi}+g^{2}_{\varphi}|||_{D(\rho-4\tau,\sigma^{2})\times U'}.
\end{equation}

Next we estimate
\begin{equation*}\label{es59}
|||s^{2}_{z}|||^{T}_{p,D(\rho-5\tau,\sigma^{2})\times U'}=\sup_{0\neq z\in l_{p}^{2}}\frac{\|\lfloor s^{2}_{z}\rceil_{D(\rho-5\tau,\sigma^{2})\times U'}(z)\|_{p}}{\|z\|_{p}},
\end{equation*}
in which
\begin{equation}\label{es60}
\begin{aligned}
\|\lfloor s^{2}_{z}\rceil_{D(\rho-5\tau,\sigma^{2})\times U'}(z)\|_{p}^{2}&=\sum_{a\in\mathcal{L}}
\left|\sum_{b\in\mathcal{L}}(\|S'^{b}_{1a}(\varphi)\|_{D(\rho-5\tau)\times U'}u_{b}+\|S'^{b}_{2a}(\varphi)\|_{D(\rho-5\tau)\times U'}v_{b})\right|^{2}\langle a\rangle^{2p}\\
&+\sum_{a\in\mathcal{L}}\left|\sum_{b\in\mathcal{L}}(\|S'^{a}_{2b}(\varphi)\|_{D(\rho-5\tau)\times U'}u_{b}+\|S'^{b}_{3a}(\varphi)\|_{D(\rho-5\tau)\times U'}v_{b})\right|^{2}\langle a\rangle^{2p}.
\end{aligned}
\end{equation}
By (\ref{es46}), we see that
$[\sum_{a\in[a]_\Delta}\sum_{b\in\mathcal{L}}\|S'^{b}_{1a}(\varphi)\|_{D(\rho-5\tau)\times U'}|u_{b}|)^{2}]^{\frac{1}{2}}$
equals to
\begin{equation*}
\left[\sum_{a\in[a]}\left(\sum_{k\in\mathbb{Z}^{\mathcal{A}}}\sum_{b\in\mathcal{L}}(|\hat{S}'^{b}_{1a}(k)|+|\partial_{w}\hat{S}'^{b}_{1a}(k)|)e^{(\rho-5\tau)|k|}|u_{b}|\right)^{2}\right]^{\frac{1}{2}},
\end{equation*}
which is less than
\begin{align*}
&\leq\sum_{k\in\mathbb{Z}^{\mathcal{A}}}\left[\sum_{a\in[a]}\left(\sum_{b\in\mathcal{L}}(|\hat{S}'^{b}_{1a}(k)|+|\partial_{w}\hat{S}'^{b}_{1a}(k)|)e^{(\rho-5\tau)|k|}|u_{b}|\right)^{2}\right]^{\frac{1}{2}}\\
&\leq\sum_{k\in\mathbb{Z}^{\mathcal{A}}}\left[\sum_{b\in\mathcal{L}}\left(\sum_{a\in[a]}(|\hat{S}'^{b}_{1a}(k)|+|\partial_{w}\hat{S}'^{b}_{1a}(k)|)^{2}\right)^{\frac{1}{2}}|u_{b}|\right]e^{(\rho-5\tau)|k|}\\
&\preceq\sum_{k\in\mathbb{Z}^{\mathcal{A}}}\frac{|k|+1}{\kappa^{2}}e^{(\rho-5\tau)|k|}\sum_{[b]\in\mathcal{E}_{\Delta}}(\|\hat{F}'^{[b]}_{1[a]}(k)+\hat{G}'^{[b]}_{1[a]}(k)\|
+\|\partial_{w}\hat{F}'^{[b]}_{1[a]}(k)+\partial_{w}\hat{G}'^{[b]}_{1[a]}(k)\|)\|u_{[b]}\| \\
&\preceq\sum_{k\in\mathbb{Z}^{\mathcal{A}}}\frac{|k|+1}{\kappa^{2}}e^{(\rho-5\tau)|k|}\sum_{a\in[a]}\sum_{b\in\mathcal{L}}(|\hat{F}'^{b}_{1a}(k)+\hat{G}'^{b}_{1a}(k)|
+|\partial_{w}\hat{F}'^{b}_{1a}(k)+\partial_{w}\hat{G}'^{b}_{1a}(k)|)\tilde{u}_{b}\\
&\preceq\frac{1}{\tau\kappa^{2}}\sum_{a\in[a]}\sum_{b\in\mathcal{L}}\|F'^{b}_{1a}(\varphi)+G'^{b}_{1a}(\varphi)\|_{D(\rho-4\tau)\times U'}\tilde{u}_{b}.
\end{align*}
In conclusion, we have
\begin{equation*}\label{es62}
\sum_{a\in[a]_\Delta}\left(\sum_{b\in\mathcal{L}}\|S'^{b}_{1a}(\varphi)\|_{D(\rho-5\tau)\times U'}|u_{b}|\right)^{2}
\preceq(\frac{d^{d}_{\Delta}}{\tau\kappa^{2}})^{2}\sum_{a\in[a]_{\Delta}}\left(\sum_{b\in\mathcal{L}}\|F'^{b}_{1a}(\varphi)+G'^{b}_{1a}(\varphi)\|_{D(\rho-4\tau)\times U'}\tilde{u}_{b}\right)^{2}.
\end{equation*}
Similar estimates hold for the other three summations in R.H.S of \eqref{es60}.
Then  we have
\begin{equation*}\label{es66}
\|\lfloor s^{2}_{z}\rceil_{D(\rho-5\tau,\sigma^{2})\times U'}(z)\|_{p}\preceq\frac{d^{d}_{\Delta}}{\tau\kappa^{2}}\|\lfloor f^{2}_{z}+g^{2}_{z}\rceil_{D(\rho-4\tau,\sigma^{2})\times U'}(\tilde{z})\|_{p}.
\end{equation*}
which together with $\|\tilde{z}\|_{p}\preceq d_{\Delta}^{d}\|z\|_{p}$ implies
\begin{equation}\label{es67}
|||s^{2}_{z}|||^{T}_{D(\rho-5\tau,\sigma^{2})\times U'}\preceq\frac{d_{\Delta}^{2d}}{\tau\kappa^{2}}|||f^{2}_{z}+g^{2}_{z}|||^{T}_{D(\rho-4\tau,\sigma^{2})\times U'}.
\end{equation}
Combining (\ref{es58}) and (\ref{es67}), we have
\begin{equation*}\label{es68}
|||X_{s^{2}}|||^{T}_{p,D(\rho-5\tau,\sigma^{2},\sigma)\times U'}\preceq \frac{d_{\Delta}^{2d}}{\tau\kappa^{2}}|||X_{f^{2}+g^{2}}|||^{T}_{p,D(\rho-4\tau,\sigma^{2},\sigma)\times U'}.
\end{equation*}

\subsubsection{\textbf{Step 3: Estimate of the vector fields $X_{s}$ and $X_{h_{1}}$.}} In this part, we shall
verify the estimates in
 (\ref{ho6}). We only estimate $X_{s}$ since $X_{h_{1}}$ is easier.

Recall that $s=s^{\varphi}+s^{0}+s^{1}+s^{2}$. By (\ref{ho1}) and (\ref{es14}), we have
\begin{equation*}\label{es69}
|||X_{s^{\varphi}}|||^{T}_{p,D(\rho-\tau,\sigma^{2},\sigma)\times U'}\preceq \frac{1}{\tau\kappa^{2}}|||X_{f^{\varphi}}|||^{T}_{p,D(\rho,\sigma^{2},\sigma)\times U'}\preceq \frac{\varepsilon}{\tau\kappa^{2}},
\end{equation*}
It then follows from Proposition \ref{p} and  (\ref{ho1}) that
\begin{equation}\label{es70}
\begin{aligned}
|||X_{\{f^{high},s^{\varphi}\}}|||&^{T}_{p,D(\rho-2\tau,(\sigma-\tau)^{2},\sigma-\tau)\times U'}\\
&\preceq
\frac{1}{\tau}|||X_{f^{high}}|||^{T}_{p,D(\rho-\tau,\sigma^{2},\sigma)\times U'}|||X_{s^{\varphi}}|||^{T}_{p,D(\rho-\tau,\sigma^{2},\sigma)\times U'}\preceq\frac{\varepsilon}{\tau^{2}\kappa^{2}},
\end{aligned}
\end{equation}
Using (\ref{ho1}), (\ref{es1}) and (\ref{es35}),
we have
\begin{equation*}\label{es71}
\begin{aligned}
&|||X_{s^{1}}|||^{T}_{p,D(\rho-3\tau,(\sigma-3\tau)^{2},\sigma-3\tau)\times U'}\preceq \frac{d_{\Delta}^{d}}{\tau\kappa^{2}}|||X_{f^{1}+g^{1}}|||^{T}_{p,D(\rho-2\tau,(\sigma-3\tau)^{2},\sigma-3\tau)\times U'}\\
\preceq &\frac{d_{\Delta}^{d}}{\tau\kappa^{2}}(|||X_{f^{1}}|||^{T}_{p,D(\rho-2\tau,(\sigma-3\tau)^{2},\sigma-3\tau)\times U'}+|||X_{g^{1}}|||^{T}_{p,D(\rho-2\tau,(\sigma-3\tau)^{2},\sigma-3\tau)\times U'})
\preceq\frac{d_{\Delta}^{d}\varepsilon}{\tau^{3}\kappa^{4}}.
\end{aligned}
\end{equation*}

Similar to that of $X_{s^{1}}$, we can estimate $X_{s^{0}}$, $X_{s^{2}}$ in sequence and finally get
\begin{equation}\label{es75}
|||X_{s}|||^{T}_{p,D(\rho-5\tau,(\sigma-5\tau)^{2},\sigma-5\tau)\times U'}\preceq\frac{d_{\Delta}^{3d}\varepsilon}{\tau^{5}\kappa^{6}}.
\end{equation}
Moreover, the vector field $X_{h_{1}}$  satisfies
\begin{equation}\label{es76}
\begin{aligned}
|||X_{h_{1}}|||&^{T}_{p,D(\rho-5\tau,(\sigma-5\tau)^{2},\sigma-5\tau)\times U'}
\leq|||X_{f^{0}+g^{0}}|||^{T}_{p,D(\rho-5\tau,(\sigma-5\tau)^{2},\sigma-5\tau)\times U'}\\
&+|||X_{f^{2}+g^{2}}|||^{T}_{p,D(\rho-5\tau,(\sigma-5\tau)^{2},\sigma-5\tau)\times U'}
\preceq\frac{d_{\Delta}^{d}\varepsilon}{\tau^{4}\kappa^{4}}.
\end{aligned}
\end{equation}

\subsubsection{\textbf{Step 4: Estimate of the vector field $X_{f_{1}}$.}}
In this step, we shall verify the estimates $\eqref{ho14}-\eqref{ho15}$.
%
Using Taylor's formula, we obtain from the homological equation (\ref{ho3}) that
\begin{align*}
&\qquad(h+f)\circ X^{t}_{s}\mid_{t=1}=(h+f^{low}+f^{high})\circ X^{t}_{s}\mid_{t=1}\\
&=h+\{h,s\}+\int_{0}^{1}(1-t)\{\{h,s\},s\}\circ X^{t}_{s}dt+f^{low}+\int_{0}^{1}\{f^{low},s\}\circ X^{t}_{s}dt \\
&+f^{high}+\{f^{high},s\}+\int_{0}^{1}(1-t)\{\{f^{high},s\},s\}\circ X^{t}_{s}dt \\
&=h+h_{1}+(1-\mathcal{T}_{\Delta'})f^{low}+f^{high}+(1-\mathcal{T}_{\Delta'})\{f^{high},s\}^{low}+\{f^{high},s\}^{high} \\
&+\int_{0}^{1}(1-t)\{\{h,s\},s\}\circ X^{t}_{s}dt+\int_{0}^{1}\{f^{low},s\}\circ X^{t}_{s}dt+\int_{0}^{1}(1-t)\{\{f^{high},s\},s\}\circ X^{t}_{s}dt.
\end{align*}
Then we have $f_{1}=f_{1}^{low}+f_{1}^{high}$, where
\begin{equation}\label{es113}
\begin{aligned}
f_{1}^{low}=&(1-\mathcal{T}_{\Delta'})f^{low}+(1-\mathcal{T}_{\Delta'})\{f^{high},s\}^{low}+\left(\int_{0}^{1}(1-t)\{\{h,s\},s\}\circ X^{t}_{s}dt\right)^{low}\\
&+\left(\int_{0}^{1}\{f^{low},s\}\circ X^{t}_{s}dt\right)^{low}+\left(\int_{0}^{1}(1-t)\{\{f^{high},s\},s\}\circ X^{t}_{s}dt\right)^{low},\\
 f_{1}^{high}=&f^{high}+\{f^{high},s\}^{high}+\left(\int_{0}^{1}(1-t)\{\{h,s\},s\}\circ X^{t}_{s}dt\right)^{high}\\
&+\left(\int_{0}^{1}\{f^{low},s\}\circ X^{t}_{s}dt\right)^{high}+\left(\int_{0}^{1}(1-t)\{\{f^{high},s\},s\}\circ X^{t}_{s}dt\right)^{high}.
\end{aligned}
\end{equation}

\emph{Firstly,} we consider the vector field generated by $(1-\mathcal{T}_{\Delta'})f^{low}$.
Observing that
\begin{equation*}\label{es78}
\begin{aligned}
\|(1-\mathcal{T}_{\Delta'})&F^{\varphi}_{\varphi}\|_{D(\rho-\tau)\times U'}=\sum_{|k|>\Delta'}(|\hat{F}^{\varphi}_{\varphi}(k;w)|+|\partial_{w}\hat{F}^{\varphi}_{\varphi}(k;w)|)e^{|k|(\rho-\tau)}\\
&\leq\sum_{|k|>\Delta'}e^{-\tau|k|}\|F^{\varphi}_{\varphi}\|_{D(\rho)\times U'}
\preceq \frac{1}{\tau^{\# \mathcal{A}}}e^{-\frac{1}{2}\tau\Delta'}\|F^{\varphi}_{\varphi}\|_{D(\rho)\times U'},
\end{aligned}
\end{equation*}
we obtain
\begin{equation*}\label{es79}
|||X_{(1-\mathcal{T}_{\Delta'})f^{\varphi}}|||^{T}_{p,D(\rho-\tau,\sigma^{2},\sigma)\times U'}
\preceq \frac{1}{\tau^{\# \mathcal{A}}}e^{-\frac{1}{2}\tau\Delta'}|||X_{f^{\varphi}}|||^{T}_{p,D(\rho,\sigma^{2},\sigma)\times U'}\preceq\frac{1}{\tau^{\# \mathcal{A}}}e^{-\frac{1}{2}\tau\Delta'}\varepsilon.
\end{equation*}
and the same estimates hold for $|||X_{(1-\mathcal{T}_{\Delta'})f^{0}}|||^{T}_{p,D(\rho-\tau,\sigma^{2},\sigma)\times U'}$
and $|||X_{(1-\mathcal{T}_{\Delta'})f^{1}}|||^{T}_{p,D(\rho-\tau,\sigma^{2},\sigma)\times U'}$.
%
%
Then we turn to
\begin{equation}\label{es83}
(1-\mathcal{T}_{\Delta'})f^{2}=(1-\mathcal{T}^{1}_{\Delta'})f^{2}+(1-\mathcal{T}^{2}_{\Delta'})f^{2},
\end{equation}
where
\begin{equation*}\label{es82}
\begin{aligned}
(1-\mathcal{T}^{1}_{\Delta'})f^{2}=\frac{1}{2}\langle \zeta,(1-\mathcal{T}_{\Delta'})F_{2}(\varphi;w)\zeta\rangle,\quad
(1-\mathcal{T}^{2}_{\Delta'})f^{2}=\frac{1}{2}\sum_{|k|>\Delta'}\langle \zeta,\mathcal{T}_{\Delta'}\hat{F}_{2}(k;w)\zeta\rangle e^{\mathrm{i}\langle k,\varphi\rangle},
\end{aligned}
\end{equation*}
It is easy to see
\begin{equation}\label{es84}
|||X_{(1-\mathcal{T}^{2}_{\Delta'})f^{2}}|||^{T}_{p,D(\rho-\tau,\sigma^{2},\sigma)\times U'}\preceq\frac{1}{\tau^{\# \mathcal{A}}}e^{-\frac{1}{2}\tau\Delta'}\varepsilon.
\end{equation}

%
%

By (\ref{ho2}), we have
\begin{equation}\label{es85}
\sup_{(\varphi,w)\in D(\rho)\times U}(|F_{2}(\varphi;w)|_{\gamma},|\partial_{w}F_{2}(\varphi;w)|_{\gamma})\leq\frac{\varepsilon}{\sigma^{2}}.
\end{equation}
Hence
\begin{equation*}\label{es86}
|\hat{F}'^{b}_{2a}(k;w)|+|\partial_{w}\hat{F}'^{b}_{2a}(k;w)|\preceq \frac{\varepsilon}{\sigma^{2}}e^{-\gamma|a-b|-\rho|k|},
\end{equation*}
which implies
\begin{equation*}\label{es87}
\begin{aligned}
\|F'^{b}_{2a\varphi}(\varphi;w)\|_{D(\rho-\tau)\times U}=&\sum_{k\in\mathbb{Z}^{\mathcal{A}}}(|\hat{F}'^{b}_{2a\varphi}(k;w)|+|\partial_{w}\hat{F}'^{b}_{2a\varphi}(k;w)|)e^{(\rho-\tau)|k|}\\
\preceq& \sum_{k\in\mathbb{Z}^{\mathcal{A}}}|k|e^{-\tau|k|}\frac{\varepsilon}{\sigma^{2}}e^{-\gamma|a-b|}
\preceq \frac{1}{\tau^{\# \mathcal{A}+1}}\frac{\varepsilon}{\sigma^{2}}e^{-\gamma|a-b|}.
\end{aligned}
\end{equation*}
Using Young's inequality (2) in \cite{EK}, we obtain
\begin{equation*}\label{es88}
|||X_{(1-\mathcal{T}^{1}_{\Delta'})f^{2}}|||^{T}_{p,D(\rho-\tau,\sigma^{2},\sigma)\times U'}\preceq\frac{1}{\gamma^{d+p}}\frac{1}{\tau^{\# \mathcal{A}+1}}\frac{\varepsilon}{\sigma^{2}}e^{-\frac{1}{2}\gamma\Delta'},
\end{equation*}
which together with \eqref{es84} leads to
\begin{equation*}\label{es89}
|||X_{(1-\mathcal{T}_{\Delta'})f^{2}}|||^{T}_{p,D(\rho-\tau,\sigma^{2},\sigma)\times U'}
\preceq\frac{1}{\tau^{\# \mathcal{A}}}e^{-\frac{1}{2}\tau\Delta'}\varepsilon+\frac{1}{\gamma^{d+p}}\frac{1}{\tau^{\# \mathcal{A}+1}}\frac{\varepsilon}{\sigma^{2}}e^{-\frac{1}{2}\gamma\Delta'}.
\end{equation*}
In conclusion, we have
\begin{equation}\label{es90}
|||X_{(1-\mathcal{T}_{\Delta'})f^{low}}|||^{T}_{p,D(\rho-\tau,\sigma^{2},\sigma)\times U'}
\preceq\frac{1}{\tau^{\# \mathcal{A}}}e^{-\frac{1}{2}\tau\Delta'}\varepsilon+\frac{1}{\gamma^{d+p}}\frac{1}{\tau^{\# \mathcal{A}+1}}\frac{\varepsilon}{\sigma^{2}}e^{-\frac{1}{2}\gamma\Delta'}.
\end{equation}

\emph{Secondly,} we consider the vector field generated by $(1-\mathcal{T}_{\Delta'})g^{low}=(1-\mathcal{T}_{\Delta'}) \{f^{high}, s\}^{low}$.
From  equation (\ref{es4}), we obtain
\begin{equation}\label{es91}
[s^{\varphi}]_{\Lambda,\gamma,\sigma;U',\rho,\mu}\preceq\frac{1}{\kappa^{2}}(\Delta')^{\exp}\varepsilon.
\end{equation}

Applying Proposition 6.6 in \cite{EK} to the equation (\ref{es5}), we obtain
\begin{equation*}\label{es92}
\|\hat{S}_1(k;\cdot)\|_{p,\gamma;U'}\preceq\frac{1}{\kappa^{2}}(\Delta\Delta')^{\exp}\|\hat{F}_1(k;\cdot)+\hat{G}_1(k;\cdot)\|_{p,\gamma;U'}.
\end{equation*}
Noticing that
\begin{equation*}\label{es93}
\{f^{high},s^{\varphi}\}=-\langle \partial_{r}f^{high},\partial_{\varphi}s^{\varphi}\rangle,
\end{equation*}
it follows from
 (\ref{ho2}), (\ref{es91}) and  in \cite[Equation (42)]{EK} that
\begin{equation}\label{es94}
[\{f^{high},s^{\varphi}\}]_{\Lambda,\gamma,\sigma;U',\rho^{(1)},\mu^{(1)}}\preceq\frac{1}{\rho-\rho^{(1)}}\frac{1}{\mu-\mu^{(1)}}\frac{1}{\kappa^{2}}(\Delta')^{\exp}\varepsilon,
\end{equation}
which together with
(\ref{ho2}) and  \eqref{es1} implies
\begin{equation*}\label{es95}
[s^{1}]_{\Lambda,\gamma,\sigma;U',\rho^{(1)},\mu^{(1)}}\preceq\frac{1}{\kappa^{4}}(\Delta\Delta')^{\exp}\frac{1}{\rho-\rho^{(1)}}\frac{1}{\mu-\mu^{(1)}}\varepsilon.
\end{equation*}
Since $s^{1}$ is independent of $r$, there is
\begin{equation}\label{es96}
[s^{1}]_{\Lambda,\gamma,\sigma;U',\rho^{(1)},\mu}\preceq\frac{1}{\kappa^{4}}(\Delta\Delta')^{\exp}\frac{1}{\rho-\rho^{(1)}}\frac{1}{\mu}\varepsilon.
\end{equation}

Next we estimate
\begin{equation*}\label{es97}
\{f^{high},s^{1}\}=-\langle \partial_{r}f^{high},\partial_{\varphi}s^{1}\rangle+\langle \partial_{\zeta}f^{high},J\partial_{\zeta}s^{1}\rangle.
\end{equation*}
By (\ref{ho2}), (\ref{es96}) and Cauchy estimates (42) in \cite{EK}, we have
\begin{equation*}\label{es98}
[\langle \partial_{r}f^{high},\partial_{\varphi}s^{1}\rangle]_{\Lambda,\gamma,\sigma;U',\rho^{(2)},\mu^{(1)}}
\preceq\frac{1}{\kappa^{4}}(\Delta\Delta')^{\exp}\frac{1}{\rho-\rho^{(1)}}\frac{1}{\rho^{(1)}-\rho^{(2)}}\frac{1}{\mu-\mu^{(1)}}\frac{1}{\mu}\varepsilon.
\end{equation*}
Applying further Proposition 3.1 (ii) in \cite{EK}, we have
\begin{equation*}\label{es99}
[\langle \partial_{\zeta}f^{high},J\partial_{\zeta}s^{1}\rangle]_{\Lambda,\gamma,\sigma^{(1)};U',\rho^{(1)},\mu}
\preceq\frac{1}{\kappa^{4}}(\Delta\Delta')^{\exp}\frac{1}{\rho-\rho^{(1)}}\frac{1}{\sigma-\sigma^{(1)}}\frac{1}{\sigma}\frac{1}{\mu}\varepsilon,
\end{equation*}
which implies
\begin{equation*}\label{es100}
[\{f^{high},s^{1}\}]_{\left\{\substack{\Lambda,\gamma,\sigma^{(1)};\\ U',\rho^{(2)},\mu^{(1)}}\right\}}
\preceq\frac{1}{\kappa^{4}}(\Delta\Delta')^{\exp}\frac{1}{\rho-\rho^{(1)}}\left(\frac{1}{\sigma-\sigma^{(1)}}\frac{1}{\sigma}+\frac{1}{\rho^{(1)}-\rho^{(2)}}\frac{1}{\mu-\mu^{(1)}}\right)\frac{1}{\mu}\varepsilon.
\end{equation*}
By (\ref{es1}) and
(\ref{es94}),
we have
\begin{equation}\label{es101}
[g^{low}]_{\left\{\substack{\Lambda,\gamma,\sigma^{(1)};\\ U',\rho^{(2)},\mu^{(1)}}\right\}}
\preceq\frac{1}{\kappa^{4}}(\Delta\Delta')^{\exp}\frac{1}{\rho-\rho^{(1)}}\left(\frac{1}{\sigma-\sigma^{(1)}}\frac{1}{\sigma}+\frac{1}{\rho^{(1)}-\rho^{(2)}}\frac{1}{\mu-\mu^{(1)}}\right)\frac{1}{\mu}\varepsilon,
\end{equation}
which implies
\begin{equation}\label{es102}
\sup_{(\varphi,w)\in D(\rho^{(2)})\times U'}\left\{
\begin{aligned} |G_{2}(\varphi;w)|_{\gamma},\\
|\partial_{w}G_{2}(\varphi;w)|_{\gamma}\end{aligned}\right\}
\preceq\frac{(\Delta\Delta')^{\exp}}{(\sigma^{(1)})^{2}\kappa^{4}}
\frac{1}{\rho-\rho^{(1)}}\left(\frac{1}{\sigma-\sigma^{(1)}}\frac{1}{\sigma}+\frac{1}{\rho^{(1)}-\rho^{(2)}}\frac{1}{\mu}\right)\frac{\varepsilon}{\mu}.
\end{equation}

By (\ref{es1})
and (\ref{es70}),
 we have
\begin{equation}\label{es103}
|||X_{g^{low}}|||^{T}_{p,D(\rho-4\tau,(\sigma-4\tau)^{2},\sigma-4\tau)\times U'}
\preceq\frac{d_{\Delta}^{d}\varepsilon}{\tau^{4}\kappa^{4}}.
\end{equation}
Following the proof of (\ref{es90}), using (\ref{es102}), (\ref{es103}), we have
\begin{equation}\label{es104}
\begin{aligned}
&|||X_{(1-\mathcal{T}_{\Delta'})g^{low}}|||^{T}_{p,D(\rho-5\tau,(\sigma-4\tau)^{2},\sigma-4\tau)\times U'}
\preceq\frac{1}{\tau^{\# \mathcal{A}}}e^{-\frac{1}{2}\tau\Delta'}\frac{d_{\Delta}^{d}\varepsilon}{\tau^{4}\kappa^{4}}\\
&\qquad+\frac{1}{\gamma^{d+p}}\frac{1}{\tau^{\# \mathcal{A}+1}}e^{-\frac{1}{2}\gamma\Delta'}\frac{(\Delta\Delta')^{\exp}}{(\sigma^{(1)})^{2}\kappa^{4}}
\frac{1}{\rho-\rho^{(1)}}\left(\frac{1}{\sigma-\sigma^{(1)}}\frac{1}{\sigma}+\frac{1}{\rho^{(1)}-\rho^{(2)}}\frac{1}{\mu}\right)\frac{\varepsilon}{\mu} \\
&\quad\preceq\frac{d_{\Delta}^{d}}{\kappa^{4}\tau^{\# \mathcal{A}+4}}e^{-\frac{1}{2}\tau\Delta'}\varepsilon
+\frac{1}{\gamma^{d+p}}\frac{1}{\tau^{\# \mathcal{A}+3}}e^{-\frac{1}{2}\gamma\Delta'}\frac{(\Delta\Delta')^{\exp}}{\sigma^{6}\kappa^{4}}\varepsilon,
\end{aligned}
\end{equation}
where $\rho^{(1)}=\rho-\tau,\rho^{(2)}=\rho-2\tau,\sigma^{(1)}=\sigma-\tau$.

\emph{Thirdly,}
by (\ref{ho1}), (\ref{ho3}), (\ref{es76}), (\ref{es103}), we have
\begin{equation}\label{es105}
|||X_{\{h,s\}}|||^{T}_{p,D(\rho-5\tau,(\sigma-5\tau)^{2},\sigma-5\tau)\times U'}\preceq\frac{d_{\Delta}^{d}\varepsilon}{\tau^{4}\kappa^{4}},
\end{equation}
which together with Proposition \ref{p} and  (\ref{es75}) implies
\begin{equation*}\label{es106}
|||X_{\{\{h,s\},s\}}|||^{T}_{p,D(\rho-6\tau,(\sigma-6\tau)^{2},\sigma-6\tau)\times U'}\preceq\frac{d_{\Delta}^{4d}\varepsilon^{2}}{\tau^{10}\kappa^{10}}.
\end{equation*}
By Proposition \ref{p}, (\ref{ho1}), (\ref{es75}), we have
\begin{equation}\label{es107}
\begin{aligned}
& |||X_{\{f^{low},s\}}|||^{T}_{p,D(\rho-6\tau,(\sigma-6\tau)^{2},\sigma-6\tau)\times U'}\preceq\frac{d_{\Delta}^{3d}\varepsilon^{2}}{\tau^{6}\kappa^{6}},\\
& |||X_{\{f^{high},s\}}|||^{T}_{p,D(\rho-6\tau,(\sigma-6\tau)^{2},\sigma-6\tau)\times U'}\preceq\frac{d_{\Delta}^{3d}\varepsilon}{\tau^{6}\kappa^{6}},
\end{aligned}
\end{equation}
which implies
\begin{equation*}\label{es109}
|||X_{\{\{f^{high},s\},s\}}|||^{T}_{p,D(\rho-7\tau,(\sigma-7\tau)^{2},\sigma-7\tau)\times U'}\preceq\frac{d_{\Delta}^{6d}\varepsilon^{2}}{\tau^{12}\kappa^{12}}.
\end{equation*}
Then applying Theorem 3.3 in \cite{CLY},
we have
\begin{equation}\label{es110}
\begin{aligned}
& |||X_{\int_{0}^{1}(1-t)\{\{h,s\},s\}\circ X^{t}_{s}dt}|||^{T}_{p,D(\rho-7\tau,(\sigma-7\tau)^{2},\sigma-7\tau)\times U'}
\preceq\frac{d_{\Delta}^{4d}\varepsilon^{2}}{\tau^{10}\kappa^{10}},\\
& |||X_{\int_{0}^{1}\{f^{low},s\}\circ X^{t}_{s}dt}|||^{T}_{p,D(\rho-7\tau,(\sigma-7\tau)^{2},\sigma-7\tau)\times U'}\preceq\frac{d_{\Delta}^{3d}\varepsilon^{2}}{\tau^{6}\kappa^{6}},\\
& |||X_{\int_{0}^{1}(1-t)\{\{f^{high},s\},s\}\circ X^{t}_{s}dt}|||^{T}_{p,D(\rho-8\tau,(\sigma-8\tau)^{2},\sigma-8\tau)\times U'}\preceq\frac{d_{\Delta}^{6d}\varepsilon^{2}}{\tau^{12}\kappa^{12}}.
\end{aligned}
\end{equation}
By (\ref{es90}), (\ref{es104}), (\ref{es110}) and (\ref{es113}), we have
\begin{equation}\label{es115}
\begin{aligned}
|||X_{f_{1}^{low}}&|||^{T}_{p,D(\rho-8\tau,(\sigma-8\tau)^{2},\sigma-8\tau)\times U'}\\
&\preceq\frac{d_{\Delta}^{d}}{\kappa^{4}\tau^{\# \mathcal{A}+4}}e^{-\frac{1}{2}\tau\Delta'}\varepsilon
+\frac{1}{\gamma^{d+p}}\frac{1}{\tau^{\# \mathcal{A}+3}}e^{-\frac{1}{2}\gamma\Delta'}\frac{(\Delta\Delta')^{\exp}}{\sigma^{6}\kappa^{4}}\varepsilon+\frac{d_{\Delta}^{6d}\varepsilon^{2}}{\tau^{12}\kappa^{12}}.
\end{aligned}
\end{equation}

\emph{Finally, }
by (\ref{ho1}), \eqref{es107},
(\ref{es110})
and \eqref{es113},
we have
\begin{equation}\label{es116}
|||X_{f_{1}^{high}}|||^{T}_{p,D(\rho-8\tau,(\sigma-8\tau)^{2},\sigma-8\tau)\times U'}\preceq 1+\frac{d_{\Delta}^{3d}\varepsilon}{\tau^{6}\kappa^{6}}+\frac{d_{\Delta}^{6d}\varepsilon^{2}}{\tau^{12}\kappa^{12}}.
\end{equation}

\subsubsection{\textbf{Step 5: Estimate of the functions $s$ and $f_{1}$}.} In this step, we shall verify the estimates \eqref{ho16}-\eqref{ho18}.

From the equation (\ref{es6}), using (\ref{ho2}) and (\ref{es101}), we obtain
\begin{equation*}\label{es117}
[s^{0}]_{\Lambda,\gamma,\sigma^{(1)};U',\rho^{(2)},\mu^{(1)}}
\preceq\frac{1}{\kappa^{6}}(\Delta\Delta')^{\exp}\frac{1}{\rho-\rho^{(1)}}\left(\frac{1}{\sigma-\sigma^{(1)}}\frac{1}{\sigma}+\frac{1}{\rho^{(1)}-\rho^{(2)}}\frac{1}{\mu-\mu^{(1)}}\right)\frac{1}{\mu}\varepsilon.
\end{equation*}
Since $s^{0}$ is independent of $\zeta$, we obtain
\begin{equation}\label{es118}
[s^{0}]_{\Lambda,\gamma,\sigma;U',\rho^{(2)},\mu^{(1)}}
\preceq\frac{1}{\kappa^{6}}(\Delta\Delta')^{\exp}\frac{1}{\rho-\rho^{(1)}}\left(\frac{1}{\sigma^{2}}+\frac{1}{\rho^{(1)}-\rho^{(2)}}\frac{1}{\mu-\mu^{(1)}}\right)\frac{1}{\mu}\varepsilon.
\end{equation}

Applying Proposition 6.7 in \cite{EK} to the equation (\ref{es7}), it follows from  (\ref{ho2}) and (\ref{es101}) that
\begin{equation}\label{es119}
\begin{aligned}
& [s^{2}]_{\left\{\substack{\Lambda'+d_\Delta+2,\gamma,\sigma^{(1)};\\
U',\rho^{(2)},\mu^{(1)}}\right\}}
\preceq\frac{1}{\kappa^{7}}(\Delta\Delta')^{\exp}\frac{1}{\rho-\rho^{(1)}}\left(\frac{1}{\sigma-\sigma^{(1)}}\frac{1}{\sigma}+\frac{1}{\rho^{(1)}-\rho^{(2)}}\frac{1}{\mu-\mu^{(1)}}\right)\frac{1}{\mu}\varepsilon,\\
& [h_{1}]_{\left\{\substack{\Lambda'+d_\Delta+2,\gamma,\sigma^{(1)};\\ U',\rho^{(2)},\mu^{(1)}}\right\}}
\preceq\frac{1}{\kappa^{4}}(\Delta\Delta')^{\exp}\frac{1}{\rho-\rho^{(1)}}\left(\frac{1}{\sigma-\sigma^{(1)}}\frac{1}{\sigma}+\frac{1}{\rho^{(1)}-\rho^{(2)}}\frac{1}{\mu-\mu^{(1)}}\right)\frac{1}{\mu}\varepsilon.
\end{aligned}
\end{equation}
Using (\ref{es91}), (\ref{es96}), (\ref{es118}) and (\ref{es119}), we obtain
\begin{equation*}\label{es121}
[s]_{\left\{\substack{\Lambda'+d_\Delta+2,\gamma,\sigma^{(1)};\\U',\rho^{(2)},\mu^{(1)}}\right\}}
\preceq\frac{1}{\kappa^{7}}(\Delta\Delta')^{\exp}\frac{1}{\rho-\rho^{(1)}}\left(\frac{1}{\sigma-\sigma^{(1)}}\frac{1}{\sigma}+\frac{1}{\rho^{(1)}-\rho^{(2)}}\frac{1}{\mu-\mu^{(1)}}\right)\frac{1}{\mu}\varepsilon.
\end{equation*}

This completes the proof of Proposition \ref{p1}.

\section{KAM theorem}\label{Sect 3}

In this section, we will prove the following KAM theorem, upon which our main result Theorem \ref{t} is an
immediate result.

\begin{thm}\label{t1}

Consider the Hamiltonian $h+f$, where
\begin{equation*}
h(r,\zeta;w)=\langle\omega(w),r\rangle+\frac{1}{2}\langle\zeta,(\Omega(w)+H(w))\zeta\rangle
\end{equation*}
satisfy (\ref{as1})-(\ref{as9}), $H(w),\partial_{w}H(w)$ are T\"{o}plitz at $\infty$ and $\mathcal{NF}_{\Delta}$ for all $w\in U$,
\begin{equation}\label{k1}
|||X_{f}|||^{T}_{p,D(\rho,\mu,\sigma)\times U}\leq\varepsilon,
\end{equation}
\begin{equation}\label{k2}
[f]_{\Lambda,\gamma,\sigma;U,\rho,\mu}\leq\varepsilon.
\end{equation}
Assume $0<\gamma,\sigma,\rho,\mu<1$, $\Lambda,\Delta\geq 3$, $\rho=\sigma$, $\mu=\sigma^{2}$, $d_{\Delta}\gamma\leq1$.
Then there is a subset $U_{\infty}\subset U$ such that if
\begin{equation*}
\varepsilon\leq c\min\left(\gamma,\rho,\frac{1}{\Delta},\frac{1}{\Lambda}\right)^{\exp},
\end{equation*}
then for all $w\in U_{\infty}$, there is a real analytic symplectic map
\begin{equation*}
\Phi:D(\frac{\rho}{2},\frac{\mu}{4},\frac{\sigma}{2})\rightarrow D(\rho,\mu,\sigma)
\end{equation*}
such that
\begin{equation*}
(h+f)\circ\Phi=h_{\infty}+f_{\infty},
\end{equation*}
where
\begin{equation*}
h_{\infty}=\langle \omega_{\infty}(w),r\rangle+\frac{1}{2}\langle\zeta,(\Omega(w)+H_{\infty}(w))\zeta\rangle,
\end{equation*}
\begin{equation*}
f_{\infty}=O(|r|^{2}+|r|\|\zeta\|_{p}+\|\zeta\|^{3}_{p})
\end{equation*}
with the estimates
\begin{align}
&\|\Phi-id\|_{p,D(\frac{\rho}{2},\frac{\mu}{4},\frac{\sigma}{2})}\leq c\varepsilon^{\frac{2}{3}},\label{k3}\\
&|||X_{f_{\infty}}|||^{T}_{p,D(\frac{\rho}{2},\frac{\mu}{4},\frac{\sigma}{2})\times U_{\infty}}
\leq c\varepsilon^{\frac{2}{3}}, \label{k4}\\
&|\omega_{\infty}(w)-\omega(w)|+|\partial_{w}(\omega_{\infty}(w)-\omega(w))|\leq c\varepsilon^{\frac{2}{3}}, \label{k5}\\
&\|H_{\infty}(w)-H(w)\|+\|\partial_{w}( H_{\infty}(w)-H(w))\|\leq c\varepsilon^{\frac{2}{3}}, \label{k6}\\
&\mathrm{meas}(U\setminus U_{\infty})\leq c \varepsilon^{\exp'},\label{k7}
\end{align}
where the exponents $\exp$, $\exp'$ depend on $d, \# \mathcal{A}, p$, the constant $c$ depends on
$d,\# \mathcal{A}, p, c_1,c_2,c_3,c_4,c_5$.

\end{thm}

Based on the KAM theorem, we prove Theorem \ref{t} on the existence and time $\delta^{-1}$ stability
of the invariant tori.
\smallskip

\noindent\textbf{Proof of Theorem \ref{t}.}
Recall the Hamiltonian formulation of NLS equation \eqref{1} in Section 1.
Let
$\omega_{a}=|a|^{2}+\hat{V}(a), a\in \mathcal{A},$
$\Omega_{a}=|a|^{2}+\hat{V}(a), a\in \mathcal{L}$,
and take $w_{a}=\hat{V}(a)$, $w\in U=[-1,1]^{\mathcal{A}}$. Then
we have
\begin{equation*}
h=\sum_{a\in \mathcal{A}}\omega_{a}r_{a}+\frac{1}{2}\sum_{a\in\mathcal{L}}\Omega_{a}(\xi_{a}^{2}+\eta_{a}^{2}),\quad
f=\varepsilon \int_{\mathbb{T}^{d}}F(x,u(x), \overline{u(x)})dx.
\end{equation*}
The T\"{o}oplitz-Lipschitz property of $f$ follows from Theorem 7.2 in \cite{EK} and the tame property follows from Section 3.5 in \cite{BG}.
By Theorem \ref{t1}, if $\varepsilon>0$ is sufficiently small, then for most $V$ (in the sense of measure), the $d$-dimensional nonlinear Schr\"odinger equation (\ref{1})
has a quasi-periodic solution.
As done in \cite{CLY}, assume $u_{0}(t,x)$ with initial value $u_{0}(0,x)$ is a quasi-periodic solution for the equation (\ref{1}),
then for any solution $u(t,x)$ with initial value $u(0,x)$ satisfying
\begin{equation*}
\|u(0,\cdot)-u_{0}(0,\cdot)\|_{H^{p}(\mathbb{T}^{d})}<\delta, \ \forall 0<\delta\ll1,
\end{equation*}
we have
\begin{equation*}
\|u(t,\cdot)-u_{0}(t,\cdot)\|_{H^{p}(\mathbb{T}^{d})}<C\delta, \ \forall 0<|t|<\delta^{-1}.
\end{equation*}
In other words, the obtained KAM tori for the nonlinear Schr\"odinger equation (\ref{1}) are of long time stability.
\qed

As remarked in Eliasson-Kuksin \cite{EK}, the size of the blocks grows  much faster than quadratically along the KAM iteration, we need to take sufficiently many
normal form computations at each step to obtain a much faster iteration scheme.

\subsection{Normal form computations.}

For $\rho_{+}<\rho$, $\gamma_{+}<\gamma$, let $\Delta'=80(\log\frac{1}{\varepsilon})^{2}\frac{1}{\min(\gamma-\gamma_{+},\rho-\rho_{+})}$, $n=[\log\frac{1}{\varepsilon}]$.
Assume $\rho=\sigma$, $\mu=\sigma^{2}$, $d_{\Delta}\gamma\leq1$. For $1\leq j\leq n$, let
\begin{equation*}
\varepsilon_{j}=\frac{\varepsilon}{\kappa^{20}}\varepsilon_{j-1}, \ \varepsilon_{0}=\varepsilon,
\end{equation*}
\begin{equation*}
\gamma_{j}=\gamma-j\frac{\gamma-\gamma_{+}}{n}, \ \gamma_{0}=\gamma,
\end{equation*}
\begin{equation*}
\rho_{j}=\rho-j\frac{\rho-\rho_{+}}{n}, \ \rho_{0}=\rho,
\end{equation*}
\begin{equation*}
\sigma_{j}=\sigma-j\frac{\sigma-\sigma_{+}}{n}, \ \sigma_{0}=\sigma,
\end{equation*}
\begin{equation*}
\mu_{j}=\sigma_{j}^{2}, \ \mu_{0}=\mu,
\end{equation*}
\begin{equation*}
\Lambda_{j}=\Lambda_{j-1}+d_{\Delta}+30, \ \Lambda_{0}=\mathrm{cte}.\max(\Lambda,d_{\Delta}^{2},d_{\Delta'}^{2}),
\end{equation*}
where the constant $\mathrm{cte}$. is the one in Proposition 6.7 in \cite{EK}.

We have the following lemma.
\begin{lem}\label{l1}
For $0\leq j<n$, consider the Hamiltonian $h+h_1+\cdots+h_{j}+f_{j}$, where
\begin{equation*}
h(r,\zeta;w)=\langle\omega(w),r\rangle+\frac{1}{2}\langle\zeta,(\Omega(w)+H(w))\zeta\rangle
\end{equation*}
satisfy \emph{(\ref{as1})-(\ref{as9})}, $H(w),\partial_{w}H(w)$ are T\"{o}plitz at $\infty$ and $\mathcal{NF}_{\Delta}$ for all $w\in U$.
Let $U' \subset U$ satisfy \emph{(\ref{sd1})-(\ref{sd4})}. For all $w\in U'$,
\begin{equation*}
h_{j}=a_{j}(w)+\langle \chi_{j}(w),r\rangle+\frac{1}{2}\langle\zeta,H_{j}(w)\zeta\rangle,
\end{equation*}
\begin{equation*}
f_{j}=f_{j}^{low}+f_{j}^{high}
\end{equation*}
satisfy
\begin{equation}\label{fi1}
|||X_{f_{j}^{low}}|||^{T}_{p,D(\rho_{j},\mu_{j},\sigma_{j})\times U'}\leq\beta^{j}\varepsilon_{j}, \  |||X_{f_{j}^{high}}|||^{T}_{p,D(\rho_{j},\mu_{j},\sigma_{j})\times U'}\leq1,
\end{equation}
\begin{equation}\label{fi2}
[f_{j}^{low}]_{\Lambda_{j},\gamma_{j},\sigma_{j};U',\rho_{j},\mu_{j}}\leq\beta^{j}\varepsilon_{j}, \ [f_{j}^{high}]_{\Lambda_{j},\gamma_{j},\sigma_{j};U',\rho_{j},\mu_{j}}\leq 1
\end{equation}
for some
\begin{equation*}
\beta\preceq \max\left(\frac{1}{\gamma-\gamma_{+}},\frac{1}{\rho-\rho_{+}},\Delta,\Lambda,\log\frac{1}{\varepsilon}\right)^{\exp_1}.
\end{equation*}
Then there exists an exponent $\exp_2$ such that if
\begin{equation*}
\varepsilon\preceq\kappa^{20}\min\left(\gamma-\gamma_{+},\rho-\rho_{+},\frac{1}{\Delta},\frac{1}{\Lambda},\frac{1}{\log\frac{1}{\varepsilon}}\right)^{\exp_2},
\end{equation*}
then for all $w\in U'$, there is a real analytic symplectic map $\Phi_{j}$ such that
\begin{equation*}
(h+h_1+\cdots+h_{j}+f_{j})\circ\Phi_{j}=h+h_1+\cdots+h_{j+1}+f_{j+1},
\end{equation*}
with the estimates
\begin{equation}\label{fi3}
|||X_{f^{low}_{j+1}}|||^{T}_{p,D(\rho_{j+1},\mu_{j+1},\sigma_{j+1})\times U'}
\preceq\beta^{j+1}\varepsilon_{j+1},
\end{equation}
\begin{equation}\label{fi4}
|||X_{f^{high}_{j+1}}|||^{T}_{p,D(\rho_{j+1},\mu_{j+1},\sigma_{j+1})\times U'}
\preceq 1+\frac{1}{\kappa^{6}}\beta^{j+1}\varepsilon_{j}+\beta^{j+1}\varepsilon_{j+1},
\end{equation}
\begin{equation}\label{fi5}
[f^{low}_{j+1}]_{\Lambda_{j+1},\gamma_{j+1},\sigma_{j+1};U',\rho_{j+1},\mu_{j+1}}
\preceq\beta^{j+1}\varepsilon_{j+1},
\end{equation}
\begin{equation}\label{fi6}
[f^{high}_{j+1}]_{\Lambda_{j+1},\gamma_{j+1},\sigma_{j+1};U',\rho_{j+1},\mu_{j+1}}
\preceq 1+\frac{1}{\kappa^{7}}\beta^{j+1}\varepsilon_{j}+\beta^{j+1}\varepsilon_{j+1},
\end{equation}
where the exponents $\exp_1$, $\exp_2$ depend on $d, \# \mathcal{A}, p$.

\end{lem}

\begin{proof}

By Proposition \ref{p1}, we can solve the homological equation
\begin{equation}\label{a1}
\{h,s_{j}\}=-\mathcal{T}_{\Delta'}f_{j}^{low}-\mathcal{T}_{\Delta'}\{f_{j}^{high},s_{j}\}^{low}+h_{j+1}
\end{equation}
with the estimates
\begin{equation}\label{a2}
\begin{aligned}
& [s_{j}]_{\left\{\substack{\Lambda_{j}+d_\Delta+2,\gamma_{j},\sigma_{j}^{(1)};\\ U',\rho_{j}^{(1)},\mu_{j}^{(1)}}\right\}}
\preceq\frac{1}{\kappa^{7}}(\Delta\Delta')^{\exp}\frac{1}{\rho_{j}-\rho_{j}^{(1)}}
\left(\frac{1}{\sigma_{j}-\sigma_{j}^{(1)}}\frac{1}{\sigma_{j}}+\frac{1}{\rho_{j}-\rho_{j}^{(1)}}\frac{1}{\mu_{j}-\mu_{j}^{(1)}}\right)\frac{1}{\mu_{j}}\beta^{j}\varepsilon_{j},\\
& [h_{j+1}]_{\left\{\substack{\Lambda_{j}+d_\Delta+2,\gamma_{j},\sigma_{j}^{(1)};\\U',\rho_{j}^{(1)},\mu_{j}^{(1)}}\right\}}
\preceq\frac{1}{\kappa^{4}}(\Delta\Delta')^{\exp}\frac{1}{\rho_{j}-\rho_{j}^{(1)}}
\left(\frac{1}{\sigma_{j}-\sigma_{j}^{(1)}}\frac{1}{\sigma_{j}}+\frac{1}{\rho_{j}-\rho_{j}^{(1)}}\frac{1}{\mu_{j}-\mu_{j}^{(1)}}\right)\frac{1}{\mu_{j}}\beta^{j}\varepsilon_{j}.
\end{aligned}
\end{equation}

Using Taylor's formula, by the homological equation (\ref{a1}), we obtain
\begin{equation}\label{a4}
(h+h_1+\cdots+h_{j}+f_{j})\circ X^{t}_{s_{j}}\mid_{t=1}=(h+h_1+\cdots+h_{j}+f_{j}^{low}+f_{j}^{high})\circ X^{t}_{s_{j}}\mid_{t=1}
\end{equation}
\begin{align*}
=&h+\{h,s_{j}\}+\int_{0}^{1}(1-t)\{\{h,s_{j}\},s_{j}\}\circ X^{t}_{s_{j}}dt+h_1+\cdots+h_{j}+\int_{0}^{1}\{h_1+\cdots+h_{j},s_{j}\}\circ X^{t}_{s_{j}}dt \\
&+f_{j}^{low}+\int_{0}^{1}\{f_{j}^{low},s_{j}\}\circ X^{t}_{s_{j}}dt+f_{j}^{high}+\{f_{j}^{high},s_{j}\}+\int_{0}^{1}(1-t)\{\{f_{j}^{high},s_{j}\},s_{j}\}\circ X^{t}_{s_{j}}dt \\
=&h+h_{1}+\cdots+h_{j+1}+(1-\mathcal{T}_{\Delta'})f_{j}^{low}+f_{j}^{high}+(1-\mathcal{T}_{\Delta'})\{f_{j}^{high},s_{j}\}^{low}+\{f_{j}^{high},s_{j}\}^{high} \\
&+\int_{0}^{1}(1-t)\{\{h,s_{j}\},s_{j}\}\circ X^{t}_{s_{j}}dt+\int_{0}^{1}\{h_1+\cdots+h_{j},s_{j}\}\circ X^{t}_{s_{j}}dt \\
&+\int_{0}^{1}\{f_{j}^{low},s_{j}\}\circ X^{t}_{s_{j}}dt+\int_{0}^{1}(1-t)\{\{f_{j}^{high},s_{j}\},s_{j}\}\circ X^{t}_{s_{j}}dt.
\end{align*}

It follows from  (\ref{fi2}) that
\begin{equation}\label{a5}
[(1-\mathcal{T}_{\Delta'})f_{j}^{low}]_{\left\{\substack{\Lambda_{j},\gamma_{j}^{(1)},\sigma_{j};\\U',\rho_{j}^{(1)},\mu_{j}}\right\}}\preceq
\left[\left(\frac{1}{\rho_{j}-\rho_{j}^{(1)}}\right)^{\#\mathcal{A}}e^{-\frac{1}{2}(\rho_{j}-\rho_{j}^{(1)})\Delta'}+e^{-(\gamma_{j}-\gamma_{j}^{(1)})\Delta'}\right]\beta^{j}\varepsilon_{j}.
\end{equation}
By (\ref{es101}), we have
\begin{equation}\label{a6}
[\{f_{j}^{high},s_{j}\}^{low}]_{\left\{\substack{\Lambda_{j},\gamma_{j},\sigma_{j}^{(1)};\\U',\rho_{j}^{(1)},\mu_{j}^{(1)}}\right\}}
\preceq\frac{1}{\kappa^{4}}(\Delta\Delta')^{\exp}\frac{1}{\rho_{j}-\rho_{j}^{(1)}}
\left(\frac{1}{\sigma_{j}-\sigma_{j}^{(1)}}\frac{1}{\sigma_{j}}+\frac{1}{\rho_{j}-\rho_{j}^{(1)}}\frac{1}{\mu_{j}-\mu_{j}^{(1)}}\right)\frac{1}{\mu_{j}}\beta^{j}\varepsilon_{j}.
\end{equation}
Hence
\begin{equation}\label{a7}
[(1-\mathcal{T}_{\Delta'})\{f_{j}^{high},s_{j}\}^{low}]_{\left\{\substack{\Lambda_{j},\gamma^{(1)}_{j},\sigma_{j}^{(1)};\\U',\rho_{j}^{(2)},\mu_{j}^{(1)}}\right\}}
\preceq\left[\left(\frac{1}{\rho_{j}^{(1)}-\rho_{j}^{(2)}}\right)^{\#\mathcal{A}}e^{-\frac{1}{2}(\rho_{j}^{(1)}-\rho_{j}^{(2)})\Delta'}+e^{-(\gamma_{j}-\gamma_{j}^{(1)})\Delta'}\right]
\end{equation}
\begin{equation*}
\times\frac{1}{\kappa^{4}}(\Delta\Delta')^{\exp}\frac{1}{\rho_{j}-\rho_{j}^{(1)}}
\left(\frac{1}{\sigma_{j}-\sigma_{j}^{(1)}}\frac{1}{\sigma_{j}}+\frac{1}{\rho_{j}-\rho_{j}^{(1)}}\frac{1}{\mu_{j}-\mu_{j}^{(1)}}\right)\frac{1}{\mu_{j}}\beta^{j}\varepsilon_{j}.
\end{equation*}

By (\ref{fi2}), (\ref{a2}), using Proposition 3.3 in \cite{EK}, we have
\begin{equation}\label{a8}
[\{f_{j}^{high},s_{j}\}]_{\Lambda_{j}+d_\Delta+5,\gamma_{j}^{(1)},\sigma_{j}^{(2)};U',\rho_{j}^{(2)},\mu_{j}^{(2)}}
\end{equation}
\begin{equation*}
\preceq\left[(\Lambda_{j}+d_\Delta+2)^{2}\left(\frac{1}{\gamma_{j}-\gamma_{j}^{(1)}}\right)^{d+1}\frac{1}{\sigma_{j}^{(1)}-\sigma_{j}^{(2)}}\frac{1}{\sigma_{j}^{(1)}}
+\frac{1}{\rho_{j}^{(1)}-\rho_{j}^{(2)}}\frac{1}{\mu_{j}^{(1)}-\mu_{j}^{(2)}}\right]
\end{equation*}
\begin{equation*}
\times\frac{1}{\kappa^{7}}(\Delta\Delta')^{\exp}\frac{1}{\rho_{j}-\rho_{j}^{(1)}}
\left(\frac{1}{\sigma_{j}-\sigma_{j}^{(1)}}\frac{1}{\sigma_{j}}+\frac{1}{\rho_{j}-\rho_{j}^{(1)}}\frac{1}{\mu_{j}-\mu_{j}^{(1)}}\right)\frac{1}{\mu_{j}}\beta^{j}\varepsilon_{j},
\end{equation*}
\begin{equation*}\label{a9}
\textrm{and}\qquad [\{f_{j}^{low},s_{j}\}]_{\Lambda_{j}+d_\Delta+5,\gamma_{j}^{(1)},\sigma_{j}^{(2)};U',\rho_{j}^{(2)},\mu_{j}^{(2)}}
\end{equation*}
\begin{equation*}
\preceq\left[(\Lambda_{j}+d_\Delta+2)^{2}\left(\frac{1}{\gamma_{j}-\gamma_{j}^{(1)}}\right)^{d+1}\frac{1}{\sigma_{j}^{(1)}-\sigma_{j}^{(2)}}\frac{1}{\sigma_{j}^{(1)}}
+\frac{1}{\rho_{j}^{(1)}-\rho_{j}^{(2)}}\frac{1}{\mu_{j}^{(1)}-\mu_{j}^{(2)}}\right]
\end{equation*}
\begin{equation*}
\times\frac{1}{\kappa^{7}}(\Delta\Delta')^{\exp}\frac{1}{\rho_{j}-\rho_{j}^{(1)}}
\left(\frac{1}{\sigma_{j}-\sigma_{j}^{(1)}}\frac{1}{\sigma_{j}}+\frac{1}{\rho_{j}-\rho_{j}^{(1)}}\frac{1}{\mu_{j}-\mu_{j}^{(1)}}\right)\frac{1}{\mu_{j}}\beta^{2j}\varepsilon_{j}^{2}.
\end{equation*}

By (\ref{a2}), (\ref{a6}), using Proposition 3.3 in \cite{EK}, we have
\begin{equation*}\label{a10}
[\{\{f_{j}^{high},s_{j}\}^{low},s_{j}\}]_{\Lambda_{j}+d_\Delta+5,\gamma_{j}^{(1)},\sigma_{j}^{(2)};U',\rho_{j}^{(2)},\mu_{j}^{(2)}}
\end{equation*}
\begin{equation*}
\begin{aligned}
\preceq& \left[(\Lambda_{j}+d_\Delta+2)^{2}\left(\frac{1}{\gamma_{j}-\gamma_{j}^{(1)}}\right)^{d+1}\frac{1}{\sigma_{j}^{(1)}-\sigma_{j}^{(2)}}\frac{1}{\sigma_{j}^{(1)}}
+\frac{1}{\rho_{j}^{(1)}-\rho_{j}^{(2)}}\frac{1}{\mu_{j}^{(1)}-\mu_{j}^{(2)}}\right]\\
&\times\frac{1}{\kappa^{11}}(\Delta\Delta')^{\exp}\frac{1}{(\rho_{j}-\rho_{j}^{(1)})^{2}}
\left(\frac{1}{\sigma_{j}-\sigma_{j}^{(1)}}\frac{1}{\sigma_{j}}+\frac{1}{\rho_{j}-\rho_{j}^{(1)}}\frac{1}{\mu_{j}-\mu_{j}^{(1)}}\right)^{2}\frac{1}{\mu_{j}^{2}}\beta^{2j}\varepsilon_{j}^{2}.
\end{aligned}
\end{equation*}

By (\ref{a2}), (\ref{a8}), using Proposition 3.3 in \cite{EK}, we have
\begin{equation*}\label{a11}
[\{\{f_{j}^{high},s_{j}\},s_{j}\}]_{\Lambda_{j}+d_\Delta+8,\gamma_{j}^{(2)},\sigma_{j}^{(3)};U',\rho_{j}^{(3)},\mu_{j}^{(3)}}
\end{equation*}
\begin{equation*}
\begin{aligned}
\preceq&\left[(\Lambda_{j}+d_\Delta+5)^{2}\left(\frac{1}{\gamma_{j}^{(1)}-\gamma_{j}^{(2)}}\right)^{d+1}\frac{1}{\sigma_{j}^{(2)}-\sigma_{j}^{(3)}}\frac{1}{\sigma_{j}^{(2)}}
+\frac{1}{\rho_{j}^{(2)}-\rho_{j}^{(3)}}\frac{1}{\mu_{j}^{(2)}-\mu_{j}^{(3)}}\right]\\
&\times\left[(\Lambda_{j}+d_\Delta+2)^{2}\left(\frac{1}{\gamma_{j}-\gamma_{j}^{(1)}}\right)^{d+1}\frac{1}{\sigma_{j}^{(1)}-\sigma_{j}^{(2)}}\frac{1}{\sigma_{j}^{(1)}}
+\frac{1}{\rho_{j}^{(1)}-\rho_{j}^{(2)}}\frac{1}{\mu_{j}^{(1)}-\mu_{j}^{(2)}}\right]\\
&\times\frac{1}{\kappa^{14}}(\Delta\Delta')^{\exp}\frac{1}{(\rho_{j}-\rho_{j}^{(1)})^{2}}
\left(\frac{1}{\sigma_{j}-\sigma_{j}^{(1)}}\frac{1}{\sigma_{j}}+\frac{1}{\rho_{j}-\rho_{j}^{(1)}}\frac{1}{\mu_{j}-\mu_{j}^{(1)}}\right)^{2}\frac{1}{\mu_{j}^{2}}\beta^{2j}\varepsilon_{j}^{2}.
\end{aligned}
\end{equation*}

By (\ref{a2})
and   \cite[Proposition 3.3]{EK}, we have $[\{h_{i+1},s_{j}\}]_{\Lambda_{j}+d_\Delta+5,\gamma_{j}^{(1)},\sigma_{j}^{(2)};U',\rho_{j}^{(2)},\mu_{j}^{(2)}}$ is less than
\begin{equation*}
\begin{aligned}
\preceq&\left[(\Lambda_{j}+d_\Delta+2)^{2}\left(\frac{1}{\gamma_{j}-\gamma_{j}^{(1)}}\right)^{d+p}\frac{\sigma_{j}^{(1)}}{\sigma_{j}^{(1)}-\sigma_{j}^{(2)}}\frac{1}{(\sigma_{i}^{(1)})^{2}}
+\frac{1}{\rho_{j}^{(1)}-\rho_{j}^{(2)}}\frac{1}{\mu_{i}^{(1)}}\right]\\
&\times\frac{1}{\kappa^{11}}(\Delta\Delta')^{\exp}\frac{1}{\rho_{j}-\rho_{j}^{(1)}}
\left(\frac{1}{\sigma_{j}-\sigma_{j}^{(1)}}\frac{1}{\sigma_{j}}+\frac{1}{\rho_{j}-\rho_{j}^{(1)}}\frac{1}{\mu_{j}-\mu_{j}^{(1)}}\right)\frac{1}{\mu_{j}}\\
&\times\frac{1}{\rho_{i}-\rho_{i}^{(1)}}\left(\frac{1}{\sigma_{i}-\sigma_{i}^{(1)}}\frac{1}{\sigma_{i}}+\frac{1}{\rho_{i}-\rho_{i}^{(1)}}\frac{1}{\mu_{i}-\mu_{i}^{(1)}}\right)\frac{1}{\mu_{i}}
\beta^{i+j}\varepsilon_{i}\varepsilon_{j}.
\end{aligned}
\end{equation*}

Take $\rho_{j}^{(l)}=\rho_{j}-\frac{l}{4}(\rho_{j}-\rho_{j+1})$, $\gamma_{j}^{(l)}=\gamma_{j}-\frac{l}{4}(\gamma_{j}-\gamma_{j+1})$,
$\sigma_{j}^{(l)}=\sigma_{j}-\frac{l}{4}(\sigma_{j}-\sigma_{j+1})$, $\mu_{j}^{(l)}=(\sigma_{j}^{(l)})^{2}$, $l=1,2,3,4$.
By (\ref{a5}), we have
\begin{equation}\label{a13}
[(1-\mathcal{T}_{\Delta'})f_{j}^{low}]_{\Lambda_{j},\gamma_{j}^{(1)},\sigma_{j};U',\rho_{j}^{(1)},\mu_{j}}\preceq
\beta\varepsilon\beta^{j}\varepsilon_{j}\preceq\beta^{j+1}\varepsilon_{j+1}.
\end{equation}
By (\ref{a7}),
we have
\begin{equation}\label{a14}
\begin{aligned}
& [(1-\mathcal{T}_{\Delta'})\{f_{j}^{high},s_{j}\}^{low}]_{\Lambda_{j},\gamma^{(1)}_{j},\sigma_{j}^{(1)};U',\rho_{j}^{(2)},\mu_{j}^{(1)}}
\preceq\frac{1}{\kappa^{4}}\beta\varepsilon\beta^{j}\varepsilon_{j}\preceq\beta^{j+1}\varepsilon_{j+1},\\
& [\{f_{j}^{high},s_{j}\}]_{\Lambda_{j}+d_\Delta+5,\gamma_{j}^{(1)},\sigma_{j}^{(2)};U',\rho_{j}^{(2)},\mu_{j}^{(2)}}
\preceq\frac{1}{\kappa^{7}}\beta\beta^{j}\varepsilon_{j}\preceq\frac{1}{\kappa^{7}}\beta^{j+1}\varepsilon_{j},\\
& [\{f_{j}^{low},s_{j}\}]_{\Lambda_{j}+d_\Delta+5,\gamma_{j}^{(1)},\sigma_{j}^{(2)};U',\rho_{j}^{(2)},\mu_{j}^{(2)}}
\preceq\frac{1}{(\Lambda_{j}+d_\Delta+5)^{14}}\beta^{j+1}\varepsilon_{j+1},\\
& [\{\{f_{j}^{high},s_{j}\}^{low},s_{j}\}]_{\Lambda_{j}+d_\Delta+5,\gamma_{j}^{(1)},\sigma_{j}^{(2)};U',\rho_{j}^{(2)},\mu_{j}^{(2)}}
\preceq\frac{1}{(\Lambda_{j}+d_\Delta+5)^{14}}\beta^{j+1}\varepsilon_{j+1},\\
& [\{\{f_{j}^{high},s_{j}\},s_{j}\}]_{\Lambda_{j}+d_\Delta+8,\gamma_{j}^{(2)},\sigma_{j}^{(3)};U',\rho_{j}^{(3)},\mu_{j}^{(3)}}
\preceq\frac{1}{(\Lambda_{j}+d_\Delta+8)^{14}}\beta^{j+1}\varepsilon_{j+1},\\
& [\{h_1+\cdots+h_{j}+(1-t)h_{j+1},s_{j}\}]_{\Lambda_{j}+d_\Delta+5,\gamma_{j}^{(1)},\sigma_{j}^{(2)};U',\rho_{j}^{(2)},\mu_{j}^{(2)}}
\preceq\frac{1}{(\Lambda_{j}+d_\Delta+5)^{14}}\beta^{j+1}\varepsilon_{j+1}.
\end{aligned}
\end{equation}

By (\ref{a1}), (\ref{a4}), we obtain
\begin{equation*}\label{a20}
\begin{aligned}
f_{j+1}=& (1-\mathcal{T}_{\Delta'})f_{j}^{low}+f_{j}^{high}+(1-\mathcal{T}_{\Delta'})\{f_{j}^{high},s_{j}\}^{low}+\{f_{j}^{high},s_{j}\}^{high}\\
&+\int_{0}^{1}(1-t)\{-\mathcal{T}_{\Delta'}f_{j}^{low},s_{j}\}\circ X^{t}_{s_{j}}dt
+\int_{0}^{1}(1-t)\{-\mathcal{T}_{\Delta'}\{f_{j}^{high},s_{j}\}^{low},s_{j}\}\circ X^{t}_{s_{j}}dt\\
&+\int_{0}^{1}\{h_1+\cdots+h_{j}+(1-t)h_{j+1},s_{j}\}\circ X^{t}_{s_{j}}dt
+\int_{0}^{1}\{f_{j}^{low},s_{j}\}\circ X^{t}_{s_{j}}dt\\
&+\int_{0}^{1}(1-t)\{\{f_{j}^{high},s_{j}\},s_{j}\}\circ X^{t}_{s_{j}}dt.
\end{aligned}
\end{equation*}
Hence
\begin{equation}\label{a21}
f^{low}_{j+1}=(1-\mathcal{T}_{\Delta'})f_{j}^{low}+(1-\mathcal{T}_{\Delta'})\{f_{j}^{high},s_{j}\}^{low}
+\left(\int_{0}^{1}(1-t)\{-\mathcal{T}_{\Delta'}f_{j}^{low},s_{j}\}\circ X^{t}_{s_{j}}dt\right)^{low}
\end{equation}
\begin{align*}
&+\left(\int_{0}^{1}(1-t)\{-\mathcal{T}_{\Delta'}\{f_{j}^{high},s_{j}\}^{low},s_{j}\}\circ X^{t}_{s_{j}}dt\right)^{low}
+\left(\int_{0}^{1}(1-t)\{\{f_{j}^{high},s_{j}\},s_{j}\}\circ X^{t}_{s_{j}}dt\right)^{low}\\
&+\left(\int_{0}^{1}\{f_{j}^{low},s_{j}\}\circ X^{t}_{s_{j}}dt\right)^{low}+\left(\int_{0}^{1}\{h_1+\cdots+h_{j}+(1-t)h_{j+1},s_{j}\}\circ X^{t}_{s_{j}}dt\right)^{low},
\end{align*}
\begin{equation}\label{a22}
\textrm{and}\qquad f^{high}_{j+1}=f_{j}^{high}+\{f_{j}^{high},s_{j}\}^{high}
+\left(\int_{0}^{1}(1-t)\{-\mathcal{T}_{\Delta'}f_{j}^{low},s_{j}\}\circ X^{t}_{s_{j}}dt\right)^{high}
\end{equation}
\begin{align*}
&+\left(\int_{0}^{1}(1-t)\{-\mathcal{T}_{\Delta'}\{f_{j}^{high},s_{j}\}^{low},s_{j}\}\circ X^{t}_{s_{j}}dt\right)^{high}
+\left(\int_{0}^{1}(1-t)\{\{f_{j}^{high},s_{j}\},s_{j}\}\circ X^{t}_{s_{j}}dt\right)^{high}\\
&+\left(\int_{0}^{1}\{f_{j}^{low},s_{j}\}\circ X^{t}_{s_{j}}dt\right)^{high}+\left(\int_{0}^{1}\{h_1+\cdots+h_{j}+(1-t)h_{j+1},s_{j}\}\circ X^{t}_{s_{j}}dt\right)^{high}.
\end{align*}

By \eqref{a14}
and \cite[Equation (53)]{EK}, we have
\begin{equation}\label{a23}
\left[\int_{0}^{1}(1-t)\{-\mathcal{T}_{\Delta'}f_{j}^{low},s_{j}\}\circ X^{t}_{s_{j}}dt\right]_{\Lambda_{j+1},\gamma_{j+1},\sigma_{j+1};U',\rho_{j+1},\mu_{j+1}}
\preceq\beta^{j+1}\varepsilon_{j+1},
\end{equation}
\begin{equation}\label{a24}
\left[\int_{0}^{1}\{f_{j}^{low},s_{j}\}\circ X^{t}_{s_{j}}dt\right]_{\Lambda_{j+1},\gamma_{j+1},\sigma_{j+1};U',\rho_{j+1},\mu_{j+1}}
\preceq\beta^{j+1}\varepsilon_{j+1},
\end{equation}
\begin{equation}\label{a25}
\left[\int_{0}^{1}(1-t)\{-\mathcal{T}_{\Delta'}\{f_{j}^{high},s_{j}\}^{low},s_{j}\}\circ X^{t}_{s_{j}}dt\right]_{\Lambda_{j+1},\gamma_{j+1},\sigma_{j+1};U',\rho_{j+1},\mu_{j+1}}
\preceq\beta^{j+1}\varepsilon_{j+1},
\end{equation}
\begin{equation}\label{a26}
\left[\int_{0}^{1}\{h_1+\cdots+h_{j}+(1-t)h_{j+1},s_{j}\}\circ X^{t}_{s_{j}}dt\right]_{\Lambda_{j+1},\gamma_{j+1},\sigma_{j+1};U',\rho_{j+1},\mu_{j+1}}
\preceq\beta^{j+1}\varepsilon_{j+1},
\end{equation}
\begin{equation}\label{a27}
\left[\int_{0}^{1}(1-t)\{\{f_{j}^{high},s_{j}\},s_{j}\}\circ X^{t}_{s_{j}}dt\right]_{\Lambda_{j+1},\gamma_{j+1},\sigma_{j+1};U',\rho_{j+1},\mu_{j+1}}
\preceq\beta^{j+1}\varepsilon_{j+1}.
\end{equation}

By (\ref{a13}), (\ref{a14}), (\ref{a21}), (\ref{a23})-(\ref{a27}), we have
\begin{equation}\label{a28}
[f^{low}_{j+1}]_{\Lambda_{j+1},\gamma_{j+1},\sigma_{j+1};U',\rho_{j+1},\mu_{j+1}}
\preceq\beta^{j+1}\varepsilon_{j+1}.
\end{equation}
By (\ref{fi2}), \eqref{a14},
(\ref{a22}), (\ref{a23})-(\ref{a27}), we have
\begin{equation*}\label{a29}
[f^{high}_{j+1}]_{\Lambda_{j+1},\gamma_{j+1},\sigma_{j+1};U',\rho_{j+1},\mu_{j+1}}
\preceq 1+\frac{1}{\kappa^{7}}\beta^{j+1}\varepsilon_{j}+\beta^{j+1}\varepsilon_{j+1}.
\end{equation*}

By (\ref{es75}), (\ref{es76}), we have
\begin{equation}\label{a30}
\begin{aligned}
& |||X_{s_{j}}|||^{T}_{p,D(\rho_{j}^{(1)},\mu_{j}^{(1)},\sigma_{j}^{(1)})\times U'}\preceq\frac{1}{\kappa^{6}}\Delta^{\exp}\left(\frac{1}{\rho_{j}-\rho_{j}^{(1)}}\right)^{5}\beta^{j}\varepsilon_{j},\\
& |||X_{h_{j+1}}|||^{T}_{p,D(\rho_{j}^{(1)},\mu_{j}^{(1)},\sigma_{j}^{(1)})\times U'}\preceq\frac{1}{\kappa^{4}}\Delta^{\exp}\left(\frac{1}{\rho_{j}-\rho_{j}^{(1)}}\right)^{4}\beta^{j}\varepsilon_{j}.
\end{aligned}
\end{equation}
By (\ref{es90}), we have
\begin{equation}\label{a32}
|||X_{(1-\mathcal{T}_{\Delta'})f_{j}^{low}}|||^{T}_{p,D(\rho_{j}^{(1)},\mu_{j}^{(1)},\sigma_{j}^{(1)})\times U'}
\end{equation}
\begin{equation*}
\preceq\left(\frac{1}{\rho_{j}-\rho_{j}^{(1)}}\right)^{\# \mathcal{A}}e^{-\frac{1}{2}(\rho_{j}-\rho_{j}^{(1)})\Delta'}\beta^{j}\varepsilon_{j}
+\frac{1}{\gamma_{j}^{d+p}}\left(\frac{1}{\rho_{j}-\rho_{j}^{(1)}}\right)^{\# \mathcal{A}+1}\frac{1}{\sigma_{j}^{2}}e^{-\frac{1}{2}\gamma_{j}\Delta'}\beta^{j}\varepsilon_{j}\preceq\beta^{j+1}\varepsilon_{j+1}.
\end{equation*}
By (\ref{es103}), (\ref{es104}), we have
\begin{equation}\label{a33}
|||X_{\{f_{j}^{high},s_{j}\}^{low}}|||^{T}_{p,D(\rho_{j}^{(1)},\mu_{j}^{(1)},\sigma_{j}^{(1)})\times U'}
\preceq\frac{1}{\kappa^{4}}\Delta^{\exp}\left(\frac{1}{\rho_{j}-\rho_{j}^{(1)}}\right)^{4}\beta^{j}\varepsilon_{j},
\end{equation}
\begin{equation}\label{a34}
\textrm{and}\qquad |||X_{(1-\mathcal{T}_{\Delta'})\{f_{j}^{high},s_{j}\}^{low}}|||^{T}_{p,D(\rho_{j}^{(1)},\mu_{j}^{(1)},\sigma_{j}^{(1)})\times U'}
\end{equation}
\begin{equation*}
\preceq\frac{1}{\kappa^{4}}(\Delta\Delta')^{\exp}\left(\frac{1}{\rho_{j}-\rho_{j}^{(1)}}\right)^{\# \mathcal{A}+4}\left[e^{-\frac{1}{2}(\rho_{j}-\rho_{j}^{(1)})\Delta'}
+\frac{1}{\gamma_{j}^{d+p}}\frac{1}{\sigma_{j}^{6}}e^{-\frac{1}{2}\gamma_{j}\Delta'}\right]\beta^{j}\varepsilon_{j}\preceq\beta^{j+1}\varepsilon_{j+1}.
\end{equation*}

By (\ref{fi1}), (\ref{a30}), using Proposition \ref{p}, we have
\begin{equation}\label{a35}
|||X_{\{f_{j}^{low},s_{j}\}}|||^{T}_{p,D(\rho_{j}^{(2)},\mu_{j}^{(2)},\sigma_{j}^{(2)})\times U'}
\preceq\frac{1}{\rho_{j}^{(1)}-\rho_{j}^{(2)}}\frac{1}{\kappa^{6}}\Delta^{\exp}\left(\frac{1}{\rho_{j}-\rho_{j}^{(1)}}\right)^{5}\beta^{2j}\varepsilon_{j}^{2},
\end{equation}
\begin{equation}\label{a36}
|||X_{\{f_{j}^{high},s_{j}\}}|||^{T}_{p,D(\rho_{j}^{(2)},\mu_{j}^{(2)},\sigma_{j}^{(2)})\times U'}
\preceq\frac{1}{\rho_{j}^{(1)}-\rho_{j}^{(2)}}\frac{1}{\kappa^{6}}\Delta^{\exp}\left(\frac{1}{\rho_{j}-\rho_{j}^{(1)}}\right)^{5}\beta^{j}\varepsilon_{j}.
\end{equation}

By (\ref{es105}), we have
\begin{equation}\label{a37}
|||X_{\{h,s_{j}\}}|||^{T}_{p,D(\rho_{j}^{(1)},\mu_{j}^{(1)},\sigma_{j}^{(1)})\times U'}
\preceq\frac{1}{\kappa^{4}}\Delta^{\exp}\left(\frac{1}{\rho_{j}-\rho_{j}^{(1)}}\right)^{4}\beta^{j}\varepsilon_{j}.
\end{equation}
By (\ref{a30}),
(\ref{a36}), (\ref{a37}), using Proposition \ref{p}, we have
\begin{equation}\label{a38}
|||X_{\{\{h,s_{j}\},s_j\}}|||^{T}_{p,D(\rho_{j}^{(2)},\mu_{j}^{(2)},\sigma_{j}^{(2)})\times U'}
\preceq\frac{1}{\rho_{j}^{(1)}-\rho_{j}^{(2)}}\frac{1}{\kappa^{10}}\Delta^{\exp}\left(\frac{1}{\rho_{j}-\rho_{j}^{(1)}}\right)^{9}\beta^{2j}\varepsilon_{j}^{2},
\end{equation}
\begin{equation}\label{a39}
|||X_{\{\{f_{j}^{high},s_{j}\},s_j\}}|||^{T}_{p,D(\rho_{j}^{(3)},\mu_{j}^{(3)},\sigma_{j}^{(3)})\times U'}
\preceq\frac{1}{\rho_{j}^{(2)}-\rho_{j}^{(3)}}\frac{1}{\rho_{j}^{(1)}-\rho_{j}^{(2)}}
\frac{1}{\kappa^{12}}\Delta^{\exp}\left(\frac{1}{\rho_{j}-\rho_{j}^{(1)}}\right)^{10}\beta^{2j}\varepsilon_{j}^{2},
\end{equation}
\begin{equation}\label{a40}
|||X_{\{h_{i+1},s_j\}}|||^{T}_{p,D(\rho_{j}^{(2)},\mu_{j}^{(2)},\sigma_{j}^{(2)})\times U'}
\preceq\frac{1}{\rho_{j}^{(1)}-\rho_{j}^{(2)}}\frac{1}{\kappa^{10}}\Delta^{\exp}\left(\frac{1}{\rho_{j}-\rho_{j}^{(1)}}\right)^{5}
\left(\frac{1}{\rho_{i}-\rho_{i}^{(1)}}\right)^{4}\beta^{i+j}\varepsilon_{i}\varepsilon_{j}.
\end{equation}

By  (\ref{a4}), we obtain
\begin{equation}\label{a41}
f_{j+1}=(1-\mathcal{T}_{\Delta'})f_{j}^{low}+f_{j}^{high}+(1-\mathcal{T}_{\Delta'})\{f_{j}^{high},s_{j}\}^{low}+\{f_{j}^{high},s_{j}\}^{high}
\end{equation}
\begin{align*}
&+\int_{0}^{1}(1-t)\{\{h,s_{j}\},s_{j}\}\circ X^{t}_{s_{j}}dt+\int_{0}^{1}\{h_1+\cdots+h_{j},s_{j}\}\circ X^{t}_{s_{j}}dt \\
&+\int_{0}^{1}\{f_{j}^{low},s_{j}\}\circ X^{t}_{s_{j}}dt+\int_{0}^{1}(1-t)\{\{f_{j}^{high},s_{j}\},s_{j}\}\circ X^{t}_{s_{j}}dt.
\end{align*}
Hence
\begin{equation}\label{a42}
\begin{aligned}
f^{low}_{j+1}=&(1-\mathcal{T}_{\Delta'})f_{j}^{low}+(1-\mathcal{T}_{\Delta'})\{f_{j}^{high},s_{j}\}^{low}+\left(\int_{0}^{1}(1-t)\{\{h,s_{j}\},s_{j}\}\circ X^{t}_{s_{j}}dt\right)^{low}\\
&+\left(\int_{0}^{1}\{h_1+\cdots+h_{j},s_{j}\}\circ X^{t}_{s_{j}}dt\right)^{low}
+\left(\int_{0}^{1}\{f_{j}^{low},s_{j}\}\circ X^{t}_{s_{j}}dt\right)^{low}\\
&+\left(\int_{0}^{1}(1-t)\{\{f_{j}^{high},s_{j}\},s_{j}\}\circ X^{t}_{s_{j}}dt\right)^{low},
\end{aligned}
\end{equation}
\begin{equation}\label{a43}
\begin{aligned}
\textrm{and}\qquad f^{high}_{j+1}=&f_{j}^{high}+\{f_{j}^{high},s_{j}\}^{high}+\left(\int_{0}^{1}(1-t)\{\{h,s_{j}\},s_{j}\}\circ X^{t}_{s_{j}}dt\right)^{high}\\
&+\left(\int_{0}^{1}\{h_1+\cdots+h_{j},s_{j}\}\circ X^{t}_{s_{j}}dt\right)^{high}
+\left(\int_{0}^{1}\{f_{j}^{low},s_{j}\}\circ X^{t}_{s_{j}}dt\right)^{high}\\
&+\left(\int_{0}^{1}(1-t)\{\{f_{j}^{high},s_{j}\},s_{j}\}\circ X^{t}_{s_{j}}dt\right)^{high}.
\end{aligned}
\end{equation}

By (\ref{a35}), (\ref{a38})-(\ref{a40}), using Theorem 3.3 in \cite{CLY}, we have
\begin{equation}\label{a44}
|||X_{\int_{0}^{1}(1-t)\{\{h,s_{j}\},s_{j}\}\circ X^{t}_{s_{j}}dt}|||^{T}_{p,D(\rho_{j+1},\mu_{j+1},\sigma_{j+1})\times U'}
\preceq\beta^{j+1}\varepsilon_{j+1},
\end{equation}
\begin{equation}\label{a45}
|||X_{\int_{0}^{1}\{h_1+\cdots+h_{j},s_{j}\}\circ X^{t}_{s_{j}}dt}|||^{T}_{p,D(\rho_{j+1},\mu_{j+1},\sigma_{j+1})\times U'}
\preceq\beta^{j+1}\varepsilon_{j+1},
\end{equation}
\begin{equation}\label{a46}
|||X_{\int_{0}^{1}\{f_{j}^{low},s_{j}\}\circ X^{t}_{s_{j}}dt}|||^{T}_{p,D(\rho_{j+1},\mu_{j+1},\sigma_{j+1})\times U'}
\preceq\beta^{j+1}\varepsilon_{j+1},
\end{equation}
\begin{equation}\label{a47}
|||X_{\int_{0}^{1}(1-t)\{\{f_{j}^{high},s_{j}\},s_{j}\}\circ X^{t}_{s_{j}}dt}|||^{T}_{p,D(\rho_{j+1},\mu_{j+1},\sigma_{j+1})\times U'}
\preceq\beta^{j+1}\varepsilon_{j+1}.
\end{equation}

By (\ref{a32}), (\ref{a34}), (\ref{a42}), (\ref{a44})-(\ref{a47}), we have
\begin{equation}\label{a48}
|||X_{f^{low}_{j+1}}|||^{T}_{p,D(\rho_{j+1},\mu_{j+1},\sigma_{j+1})\times U'}
\preceq\beta^{j+1}\varepsilon_{j+1}.
\end{equation}
By (\ref{fi1}), (\ref{a36}), (\ref{a43}), (\ref{a44})-(\ref{a47}), we have
\begin{equation}\label{a49}
|||X_{f^{high}_{j+1}}|||^{T}_{p,D(\rho_{j+1},\mu_{j+1},\sigma_{j+1})\times U'}
\preceq 1+\frac{1}{\kappa^{6}}\beta^{j+1}\varepsilon_{j}+\beta^{j+1}\varepsilon_{j+1}.
\end{equation}

\end{proof}

\subsection{Iterative Lemma}

Assume $\rho=\sigma$, $\mu=\sigma^{2}$, $d_{\Delta}\gamma\leq1$. For $m\geq 0$, let
\begin{equation*}
\varepsilon_{m}=e^{-\frac{1}{20}(\log\frac{1}{\varepsilon_{m-1}})^{2}}, \ \varepsilon_{0}=\varepsilon,
\end{equation*}
\begin{equation*}
\vartheta_{m}=\frac{\sum_{j=1}^{m}\frac{1}{j^{2}}}{2\sum_{j=1}^{\infty}\frac{1}{j^{2}}}, \ \vartheta_{0}=0,
\end{equation*}
\begin{equation*}
\rho_{m}=(1-\vartheta_{m})\rho, \ \rho_{0}=\rho,
\end{equation*}
\begin{equation*}
\sigma_{m}=(1-\vartheta_{m})\sigma, \ \sigma_{0}=\sigma,
\end{equation*}
\begin{equation*}
\mu_{m}=\sigma_{m}^{2}, \ \mu_{0}=\mu,
\end{equation*}
\begin{equation*}
\gamma_{m}=d^{-1}_{\Delta_{m}}, \ \gamma_{0}=\min(\gamma,d^{-1}_{\Delta}),
\end{equation*}
\begin{equation*}
\Delta_{m}=80(\log\frac{1}{\varepsilon_{m-1}})^{2}\frac{1}{\min(\gamma_{m-1},\rho_{m-1}-\rho_{m})}, \ \Delta_{0}=\Delta,
\end{equation*}
\begin{equation*}
\Lambda_{m}=\mathrm{cte}.d^{2}_{\Delta_{m}},
\end{equation*}
where the constant $\mathrm{cte}$. is the one in Proposition 6.7 in \cite{EK}.

We have the following iterative lemma.

\begin{lem}\label{l2}
For $m\geq0$, consider the Hamiltonian $h_{m}+f_{m}$, where
\begin{equation*}
h_{m}=\langle \omega_{m}(w),r\rangle+\frac{1}{2}\langle\zeta,(\Omega(w)+H_{m}(w))\zeta\rangle,
\end{equation*}
$H_{m}(w),\partial_{w}H_{m}(w)$ are T\"{o}plitz at $\infty$ and $\mathcal{NF}_{\Delta_{m}}$ for all $w\in U_{m}$,
\begin{equation*}
f_{m}=f_{m}^{low}+f_{m}^{high}
\end{equation*}
satisfy
\begin{equation}\label{ii1}
|||X_{f_{m}^{low}}|||^{T}_{p,D(\rho_{m},\mu_{m},\sigma_{m})\times U_{m}}\leq\varepsilon_{m},
 \  |||X_{f_{m}^{high}}|||^{T}_{p,D(\rho_{m},\mu_{m},\sigma_{m})\times U_{m}}\leq \varepsilon+\sum_{j=1}^{m}\varepsilon_{j-1}^{\frac{2}{3}},
\end{equation}
\begin{equation}\label{ii2}
[f_{m}^{low}]_{\Lambda_{m},\gamma_{m},\sigma_{m};U_{m},\rho_{m},\mu_{m}}\leq\varepsilon_{m},
\ [f_{m}^{high}]_{\Lambda_{m},\gamma_{m},\sigma_{m};U_{m},\rho_{m},\mu_{m}}\leq \varepsilon+\sum_{j=1}^{m}\varepsilon_{j-1}^{\frac{2}{3}}.
\end{equation}
Assume for all $w\in U_{m}$,
\begin{equation}\label{ii3}
\begin{aligned}
& |\omega_{m}(w)-\omega(w)|+|\partial_{w}(\omega_{m}(w)-\omega(w))|\leq\sum_{j=1}^{m}\varepsilon_{j-1}^{\frac{2}{3}},\\
&\|H_{m}-H\|_{U_{m}}+\langle H_{m}-H\rangle_{\Lambda_{m};U_{m}}\leq\sum_{j=1}^{m}\varepsilon_{j-1}^{\frac{2}{3}}.
\end{aligned}
\end{equation}
Then there is a subset $U_{m+1}\subset U_{m}$ such that if
\begin{equation*}
\varepsilon\preceq\min\left(\gamma,\rho,\frac{1}{\Delta},\frac{1}{\Lambda}\right)^{\exp},
\end{equation*}
then for all $w\in U_{m+1}$, there is a real analytic symplectic map $\Phi_{m}$ such that
\begin{equation*}
(h_{m}+f_{m})\circ\Phi_{m}=h_{m+1}+f_{m+1}
\end{equation*}
with the following estimates
\begin{align*}
& |||X_{f^{low}_{m+1}}|||^{T}_{p,D(\rho_{m+1},\mu_{m+1},\sigma_{m+1})\times U_{m+1}}
\leq\varepsilon_{m+1},\\
& |||X_{f^{high}_{m+1}}|||^{T}_{p,D(\rho_{m+1},\mu_{m+1},\sigma_{m+1})\times U_{m+1}}
\leq \varepsilon+\sum_{j=1}^{m+1}\varepsilon_{j-1}^{\frac{2}{3}},\\
& [f^{low}_{m+1}]_{\Lambda_{m+1},\gamma_{m+1},\sigma_{m+1};U_{m+1},\rho_{m+1},\mu_{m+1}}
\leq\varepsilon_{m+1},\\
& [f^{high}_{m+1}]_{\Lambda_{m+1},\gamma_{m+1},\sigma_{m+1};U_{m+1},\rho_{m+1},\mu_{m+1}}
\leq \varepsilon+\sum_{j=1}^{m+1}\varepsilon_{j-1}^{\frac{2}{3}},\\
& |\omega_{m+1}(w)-\omega(w)|+|\partial_{w}(\omega_{m+1}(w)-\omega(w))|\leq\sum_{j=1}^{m+1}\varepsilon_{j-1}^{\frac{2}{3}},\\
&\|H_{m+1}-H\|_{U_{m+1}}+\langle H_{m+1}-H\rangle_{\Lambda_{m+1};U_{m+1}}\leq\sum_{j=1}^{m+1}\varepsilon_{j-1}^{\frac{2}{3}},\\
&\mathrm{meas}(U_{m}\setminus U_{m+1})\preceq \varepsilon_{m}^{\exp'},
\end{align*}
where the exponents $\exp$, $\exp'$ depend on $d, \# \mathcal{A}, p$.

\end{lem}

\begin{proof}
Take $\kappa^{20}=\varepsilon^{\frac{1}{20}}$ in Lemma \ref{l1},
there is a real analytic symplectic map $\Phi$ such that
\begin{equation*}
(h+f)\circ\Phi=h+h_1+\cdots+h_{n}+f_{n}.
\end{equation*}
Let $h_{+}=h+h_1+\cdots+h_{n}$, $f_{+}=f_{n}$.
Using Lemma \ref{l1}, we prove the iterative lemma.
\end{proof}

Now Theorem \ref{t1} follows from Lemma \ref{l2} and we omit its proof here.

\end{document}